\documentclass[11pt]{amsart}

\usepackage[margin=.8in, top=1.1in, bottom=0.9in]{geometry}
\usepackage{amssymb}
\usepackage{amsmath,amscd}
\usepackage{mathtools}
\usepackage{mathrsfs}
\usepackage[hidelinks]{hyperref}
\usepackage{url}
\usepackage{enumerate}
\usepackage{caption}
\usepackage{tikz}
\usepackage{tikz-cd}
\usepackage{pdfpages}
\usepackage{color, colortbl}
\usepackage{array}
\usepackage{tabularray}
\usepackage{stmaryrd}
\usepackage{bm}
\usepackage{varioref}
\usepackage{cleveref}
\usepackage{xcolor}

\usepackage{wrapfig}
\usepackage{cancel}
\usepackage{blkarray}
\newcommand{\matindex}[1]{\mbox{\scriptsize#1}}

\newtheorem{thm}{Theorem}

\theoremstyle{definition}
\newtheorem{definition}[thm]{Definition}

\newtheorem{Fact}[thm]{Fact}

\newtheorem{Theorem}[thm]{Theorem}
\newenvironment{theorem}
  {\begin{Theorem}}{\end{Theorem}}
  
\newtheorem{Lemma}[thm]{Lemma}
\newenvironment{lemma}
  {\begin{Lemma}}{\end{Lemma}}

\newtheorem{Remark}[thm]{Remark}

\newtheorem{Proposition}[thm]{Proposition}
\newenvironment{proposition}
  {\begin{Proposition}}{\end{Proposition}}

\newtheorem{Corollary}[thm]{Corollary}

\newtheorem{Conjecture}[thm]{Conjecture}

\crefname{theorem}{theorem}{theorems}
\Crefname{theorem}{Theorem}{Theorems}
\crefname{lemma}{lemma}{lemmas}
\Crefname{lemma}{Lemma}{Lemmas}
\crefname{proposition}{proposition}{propositions}
\Crefname{proposition}{Proposition}{Propositions}
\crefname{example}{example}{examples}
\Crefname{example}{Example}{Examples}
\crefname{diagram}{diagram}{diagrams}
\Crefname{diagram}{Diagram}{Diagrams}

\newcommand{\A}{\mathcal{A}}

\newcommand{\E}{\mathcal{E}}
\newcommand{\K}{\mathcal{K}}
\newcommand{\U}{\mathcal{U}}

\newcommand{\V}{\mathcal{V}}
\newcommand{\J}{\mathcal{J}}

\newcommand{\B}{\mathcal{B}}

\newcommand{\I}{\mathcal{I}}
\newcommand{\C}{\mathcal{C}}

\newcommand{\F}{\mathcal{F}}

\newcommand{\NN}{\mathbb{N}}
\newcommand{\ZZ}{\mathbb{Z}}

\newcommand{\RR}{\mathbb{R}}

\DeclareMathOperator{\HHH}{H}

\DeclareMathOperator{\id}{id}

\DeclareMathOperator{\sgn}{sgn}

\DeclareMathOperator{\Hom}{Hom}
\newcommand{\HomZ}[1]{\Hom(#1,\ZZ)}
\DeclareMathOperator{\Alt}{Alt}

\DeclarePairedDelimiter{\squareBr}{ [ }{ ] }

\newcommand{\mapstoUP}[2]{
    \begin{array}{c} #2 \\[-.13cm] \rotatebox{90}{$\mapsto$} \\[-.33cm] #1 \end{array}
}
\newcommand{\mapstoUPleft}[2]{
    \begin{array}{l} #2 \\[-.13cm] \hspace{.4cm}\rotatebox{90}{$\mapsto$} \\[-.33cm] #1 \end{array}
}

\newcommand{\uZZ}{\underline{\mathbb{Z}}}
\newcommand{\uRR}{\underline{\mathbb{R}}}

\newcommand{\h}[1]{\widehat{#1}}
\renewcommand{\tilde}{\widetilde}
\newcommand{\oline}[1]{\overline{#1}}

\newcommand{\bmi}{\bm{i}} 
\newcommand{\bmj}{\bm{j}}

\newcommand{\Cech}{\v{C}ech}
\newcommand{\cdelta}{\check{\delta}}
\newcommand{\tdelta}{\tilde{\delta}}
\newcommand{\cpartial}{\check{\partial}}

\newcommand{\sCC}{C^{\mathrm{sin}}}
\newcommand{\CCs}{C_{\mathrm{sin}}}

\newcommand{\cHH}{\check{\mathrm{H}}}
\newcommand{\dHH}{\mathrm{H}^{\Delta}}
\newcommand{\HHd}{\mathrm{H}_{\Delta}}
\newcommand{\sHH}{\mathrm{H}^{\mathrm{sin}}}
\newcommand{\HHs}{\mathrm{H}_{\mathrm{sin}}}
\newcommand{\HHdR}{\mathrm{H}_{\mathrm{dR}}}

\newcommand{\HH}{\mathrm{H}}

\newcommand{\Cc}{\check{C}}
\newcommand{\dC}{C^{\Delta}}
\newcommand{\Cd}{C_{\Delta}}

\newcommand{\DS}{\displaystyle}
\newcommand{\PO}{\displaystyle\prod}
\newcommand{\POO}[1]{\displaystyle\prod_{|\bmi|=#1}}
\newcommand{\POOj}[1]{\displaystyle\prod_{|\bmj|=#1}}
\newcommand{\BO}{\displaystyle\bigoplus}
\newcommand{\BOO}[1]{\displaystyle\bigoplus_{|\bmi|=#1}}

\newcommand{\tCs}{\widetilde{C}_{sin}}
\newcommand{\Cs}{C_{sin}}
\newcommand{\sC}{C^{sin}}

\title{The evaluation isomorphism of singular cohomology on the \v{C}ech nerve}

\author{Marco Belli}

\begin{document}

\begin{abstract}
    We show that the canonical isomorphism $\HHs^*(X,\uZZ)\xrightarrow{\sim} \cHH_{\U}^*(X,\uZZ)$  is the evaluation of singular cohomology classes at the simplices of the \Cech{} nerve $N\U$ through a homotopy equivalence $N\U\to X$. 

    We apply this to real tori $V/\Lambda$ to show that the automorphy map $\HHs^k(V/\Lambda, \ZZ) \xrightarrow{\sim} \HH^k(\Lambda, \ZZ)$ and the group cohomology map $\HH^k(\Lambda, \ZZ) \xrightarrow{\sim} \cHH^k_{\U}(V/\Lambda,\uZZ)$ induced by the universal cover are inverse to each other.
    
    {\small\[\hspace{-.1in} \begin{tikzcd}
        \HHs^k(X,\uZZ)
        \arrow[r,"\sim"',"\substack{\text{sheaf}\\\text{cohomology}}"]
        &[.7cm] \cHH^k_{\U}(X,\uZZ) &[1.2cm]
        \HHs^k(V/\Lambda,\uZZ)
        \arrow[r,"\sim"',"\substack{\text{sheaf}\\\text{cohomology}}"]
        &[.8cm] \cHH^k_{\U}(V/\Lambda,\uZZ)
        \\[.2cm]
        \HHs^k(X,\ZZ) \arrow[u,"\sim"{sloped},"\text{sheafification}"'] \arrow[r,"\sim"'{sloped},"\text{restriction}"] & \HH_{\Delta}^k(N\U,\ZZ) 
        \arrow[u,"\sim"{sloped},"\mapstoUP{\squareBr{\alpha}}{ (-1)^k\squareBr{(\alpha(x_{\bmi_0}\cdots x_{\bmi_k}))_{\bmi}} }"'] &
        \HHs^k(V/\Lambda,\ZZ) \arrow[u,"\sim"{sloped},"\text{sheafification}"'] \arrow[r,"\sim"'{sloped}," \text{automorphy}"] & 
        \HHH^k(\Lambda,\ZZ)  
        \arrow[u,"\sim"{sloped},"\mapstoUP{\squareBr{F}}{ \squareBr{ (F(\lambda_{\bmi_0\bmi_1},\dots, \lambda_{\bmi_{k-1},\bmi_k}))_{\bmi}} }"'] 
    \end{tikzcd}\] \normalsize}
\end{abstract}

\maketitle


\section{\large Introduction}

Let $X$ be a semi-locally contractible topological space with a good open cover $\U$.
\footnote{With semi-locally contractible we understand that any open set $U\subseteq X$ has an open cover $\{W_i\}_i$ such that the inclusions $W_i\hookrightarrow U$ are homotopic to a constant map; see \cite{sella}. With a good open cover we understand an open cover $\{U_i\}_i$ such that every non-empty finite intersection of the $U_i$ is a contractible topological space.}
There are two fundamental acyclic resolutions of the constant sheaf $\uZZ$ on $X$, namely  by the sheaves of singular cochains $\tilde{C}_{sin}^*$ and of \Cech{} cochains $\check{C}_{\U}^*$, used for their geometrical and combinatorial versatility respectively.
\footnote{
The complex of sheaves $\tilde{C}_{sin}^*$ is obtained by sheafifying the presheaf of singular cochains $C_{sin}^*(\bullet,\ZZ)$, and semi-local contractability implies that it is an acyclic resolution of $\uZZ$.
With $\HHs^k(X,\ZZ)$ we denote the classical definition of singular cohomology \cite{hatcher}, whereas with $\HHs^k(X,\uZZ)$ we denote the sheaf cohomology of $\uZZ$ computed with the acyclic resolution $\tilde{C}_{sin}^*$.
It is a classical result that, if $X$ is semi-locally contractible and paracompact, then the sheafification map $C_{sin}^*(X)\to \tilde{C}_{sin}^*(X)$ is a chain homotopy equivalence, so $\HHs^k(X,\ZZ)\cong \HHs^k(X,\uZZ)$.
However, the paracompactness assumption is superflous, as shown in \cite{sella} and \cite{petersen} by two different methods.
}

The two cohomologies relate to the ones computed by injective resolutions $\HH^*(X,\uZZ)$, unique up to unique isomorphism, through natural isomorphisms 
$ \HHs^*(X,\uZZ)\xrightarrow{\sim} \HH^*(X,\uZZ) \xleftarrow{\sim} \cHH_{\U}^*(X,\uZZ) $
defined from the universal $\delta$-functor properties.
The drawback is that the isomorphism $\HHs^*(X,\uZZ)\xrightarrow{\sim} \cHH_{\U}^*(X,\uZZ)$ thus induced is not computable. So one falls back to relating singular cohomology to a more concrete variant, from which the combinatorial information of \Cech{} cohomology can be extracted more explicitely.

For example, if $X$ is a smooth manifold, the isomorphism $\cHH_{\U}^*(X,\uRR)\xrightarrow{\sim}\HH_{dR}^*(X,\uRR)$ is described explicitely by the collating formula \cite[Proposition 9.8]{botttu}.
Moreover, there is the chain homotopy between the acyclic resolutions $C_{dR}^*$ and $\tilde{C}^*_{sin}$, given by integration of differential forms, which describes the isomorphism $\HH_{dR}^*(X,\uRR)\xrightarrow{\sim}\HHs^*(X,\uRR)$. In short, we have the following commutative diagram, where the maps to $\HH^*(X,\uRR)$ are the ones induced by sheaf cohomology.

\[\begin{tikzcd}[column sep=large]
    \HHs^*(X,\uRR) \arrow[r,"\sim"] & \HH^*(X,\uRR) &
    \cHH^*_{\U}(X,\uRR) \arrow[l,"\sim"'] \arrow[dl, "\sim"{sloped}, "\substack{\text{collating}\\\text{formula}}"] \\[.2cm]
    \HHs^*(X,\RR) \arrow[u,"\sim"{sloped},"\text{sheafification}"'] &
    \HH_{dR}^*(X,\uRR) \arrow[u,"\sim"{sloped}] \arrow[l,"\sim"'{sloped},"\int_{\bullet}\omega\,\mapsfrom\,\omega"']
\end{tikzcd}\]


We explore an analogous scenario for integral coefficients.
\Cech{} cohomology is, almost by definition, the $\Delta$-complex cohomology of the simplicial complex given by the nerve $N\U$ of the good cover.
\footnote{We follow \cite{hatcher} for our $\Delta$-complex cohomology theory.}
On the other hand, for any good cover $\U$ there exists a homotopy equivalence $\iota:N\U\to X$, at least if $X$ is paracompact \cite[Corollary 4G.3]{hatcher}. It is thus natural to ask what is the relationship of these two characterisations with the isomorphism induced by sheaf cohomology. First of all, $\iota$ can't be just any homotopy equivalence, it has to map the simplicial complex $N\U$ in $X$ in a way consistent with the cover $\U$. The expected result that we will prove is that they are related by the identity, up to a sign $(-1)^k$.
\begin{equation} \label{eq:MainDiagram} \begin{tikzcd}
    \HHs^*(X,\uZZ) \arrow[r,"\sim"]
    & \HH^*(X,\uZZ) & \arrow[l,"\sim"'] \cHH^*_{\U}(X,\uZZ) & \hspace{1.2cm}&
    \HH_{\Delta}^*(N\U,\ZZ) \arrow[ll,"\sim"," {[ \,(\alpha((x_{\bmi_0}\cdots x_{\bmi_k})))_{\bmi_0,\dots,\bmi_k} ]} \mapsfrom {[\alpha]} "'] \\
    \HHs^*(X,\ZZ) \arrow[u,"\sim"{sloped},"\text{sheafification}"'] \arrow[rr,"\sim"'{sloped}," \iota^*"]
    && \HHs^*(N\U,\ZZ) \arrow[rr,"\sim"',"\text{restriction}"] &&
    \HH_{\Delta}^*(N\U,\ZZ) \arrow[u,dotted,no head,"?"']
\end{tikzcd}\end{equation}

Let $X$ be a topological space with an open cover $\U=\{U_i\}_{i}$.
We will denote in bold $\bmi$ non-empty finite sets of distinct indices $\{\bmi_0,...,\bmi_k\}$,  $|\bmi|:=k$, of the cover $\U$ such that the intersection $U_{\bmi}:=\bigcap_{i\in\bmi} U_i$ is non-empty.
Denote  by $x_i$ the $0$-simplex of $N\U$ corresponding to $U_i$. If the cover $\U$ is ordered, the convention we adopt for the simplicial structure on $N\U$ is that every $1$-simplex corresponding to  $U_i\cap U_j$ with $i<j$ is oriented from $x_i$ to $x_j$, we thus denote such a $1$-simplex by $(x_i\,x_j):\Delta^1\hookrightarrow N\U$.
It then follows that any $m$-simplex of $N\U$ corresponding to a non-empty finite intersection $U_{\bmi}$ with $\bmi_0<...<\bmi_m$ has vertices in order $x_{\bmi_0},...,x_{\bmi_m}$, we thus denote such an $m$-simplex by $(x_{\bmi_0} \cdots x_{\bmi_{m}}):\Delta^m\hookrightarrow N\U$. The boundary operator then takes on the expression
$$ \partial (x_{\bmi_0} \cdots x_{\bmi_{m}}) = \sum_{j=0}^{|\bmi|} (-1)^j (x_{\bmi_0} \cdots \widehat{x_{\bmi_{j}}} \cdots x_{\bmi_{m}}) .$$

\begin{definition}
    We say that a homotopy equivalence  $\iota:N\U\to X$ \emph{respects $\U$} if every simplex $(x_{\bmi_0} \cdots x_{\bmi_{m}})$ of $N\U$ has image $\iota\circ(x_{\bmi_0} \cdots x_{\bmi_{m}})$ contained in $U_{\bmi_0}\cup\cdots\cup U_{\bmi_m}$.
\end{definition}

The next lemma is what allows one to apply Theorem \ref{thm:Theorem} in general.
\begin{lemma} \label[lemma]{thm:HomotopyExistence}
    If $X$ is paracompact and $\U$ is a good open cover, then there exists a homotopy equivalence $\iota:N\U\to X$ that respects $\U$.
\end{lemma}
The existence of a homotopy equivalence is the classic nerve theorem \cite[Corollary 4G.3]{hatcher}, but for an explicit construction that respects the cover $\U$ we refer to \cite[Theorem 2.1.(1)]{paris}.
    
\begin{theorem} \label{thm:Theorem}
    Let $X$ be semi-locally contractible and paracompact, 
    \footnote{ The paracompactness assumption is only necessary to potentially extend the homotopy $\iota$ to a bigger \Cech{} nerve using \Cref{thm:HomotopyExistence}. }
    $\U$ an ordered good open cover of $X$, and $\iota:N\U\to X$ a homotopy equivalence that respects $\U$. 
    Then, the isomorphism $\HHs^*(X,\ZZ)\xrightarrow{\sim}\cHH_{\U}^*(X,\uZZ)$ induced by sheaf cohomology precomposed with sheafification
    is given by
    $$ [\alpha] \in \HHs^k(X,\ZZ) \;\mapsto\;
    (-1)^k[ \,(\alpha(\iota\circ (x_{\bmi_0}\cdots x_{\bmi_k})))_{\bmi_0,...,\bmi_k} ] \in\cHH_{\U}^k(X,\uZZ),\quad \text{ for }k\in\NN,    $$
    where $(x_{\bmi_0}\cdots x_{\bmi_k})$ is the $k$-simplex of the simplicial complex $N\U$ corresponding to $U_{\bmi_0}\cap \cdots \cap U_{\bmi_k}$. In other words, the isomorphism missing in diagram \eqref{eq:MainDiagram} is $(-1)^k$ the identity.
\end{theorem}



\section{\large Prelude on Homological Algebra} \label{HomologicalSection}

The fact that sheaf cohomology can be computed using acyclic resolutions is a standard result in homological algebra, but the relation between multiple acyclic resolutions is not usually explained; for completeness, we discuss it here, requiring us to unravel some well-known definitions.

Let $X$ be a topological space and $\F$ a sheaf of abelian groups on $X$.
Sheaf cohomology is, up to unique isomorphism, given by
$ \HH^k(X,\F) \coloneq \HH^k(\I^*(X)), $
where $0\to\F\to \I^* $ is an injective resolution of $\F$. This definition gives a universal $\delta$-functor extending the global sections functor. 

If we are given a merely acyclic resolution $0\to\F\to \A^* $, by breaking it up in SESs
and using the LES in cohomology inductively one obtains the isomorphisms
$ \HH^k(\A^*(X)) \xrightarrow{\sim}  \HH^k(\I^*(X)) $;
see for example \cite[Remark 10.2]{wedhorn}.
An equivalent way of obtaining these isomorphisms is the following. Take an injective resolution of the whole complex of sheaves $0\to\F\to \A^* $, meaning a square diagram $\J^{*,*}$ with horizontal and vertical differentials $d_{\I}$, $d_{\A}$ and injective resolutions $\I^a \hookrightarrow \J^{a,*}$, $\A^b \hookrightarrow \J^{*,b}$ for every $a,b$. 
Taking global sections, one obtains a diagram with exact rows and columns, except the first row and first column: 
\begin{equation} \label{eq:AcyclicCoh} \begin{tikzcd} [column sep=small, row sep=small]
    \A^* \arrow[r,hook] & \J^{*,*} &{}&  
    \A^*(X) \arrow[r,hook] & \J^{*,*}(X) \\
    \F \arrow[r,hook]\arrow[u,hook] & \I^* \arrow[u,hook]
    &{}& 
    \F(X) \arrow[r,hook]\arrow[u,hook] & \I^*(X) \arrow[u,hook] \\
\end{tikzcd}\end{equation}
Consider the double complex $J^*=\bigoplus_{a+b=*}\J^{a,b}(X)$ with differential acting on $\J^{a,b}(X)$ by 
$d_{\I} + (-1)^a d_{\A}$. 
Then, by the exactness of the rows and columns, the inclusions $\A^*(X)\hookrightarrow J^* \hookleftarrow \I^*(X)$ are chain maps which induce isomorphisms
$$ \HH^k(\A^*(X)) \xrightarrow{\sim} \HH^k(J^*) \xleftarrow{\sim} \HH^k(\I^*(X)) \eqcolon \HH^k(X,\F). $$
Concretely, the induced isomorphism $\HH^k(\A^*(X))\xrightarrow{\sim}\HH^k(\I^*(X))$ is obtained by chasing cocycles across the anti-diagonals of diagram \eqref{eq:AcyclicCoh}.

Let now $0\to\F\to \B^* $ be a second acyclic resolution, with a chain map $\phi:\A^*\to\B^*$ which restricts to the identity on the subsheaf $\F$ and induces a quasi-isomorphism $\phi(X):\A^*(X)\xrightarrow{\text{q.i.s.}}\B^*(X)$ on global sections. Similarly as above, we can take an injective resolution of both complexes of sheaves, with the resolution $0\to\F\to \I^* $ common to both. The chain map $\phi$ extends to a chain map of the square diagrams:
\begin{equation*} \label{eq:AcyclicMor} \begin{tikzcd} [column sep=large, row sep=small, outer sep=-2.2pt]
    & \B^* \arrow[r,hook] & \K^{*,*}
    {}& 
    & \B^*(X)\arrow[r,hook] & \K^{*,*}(X)  \\
    \A^* \arrow[r,hook]\arrow[ur,pos=0.4,"\phi"] & \J^{*,*} \arrow[ur,pos=0.4,"\phi"]
    &{}&  
    \A^*(X) \arrow[r,hook]\arrow[ur,pos=0.4,"\phi(X)"] & \J^{*,*}(X) \arrow[ur,pos=0.4,"\phi(X)"] \\
    \F \arrow[r,hook] \arrow[u,hook] \arrow[uur, hook, bend right=6, crossing over] & 
    \I^* \arrow[u,hook] \arrow[uur, hook, bend right=6, crossing over]
    &{}&
    \F(X) \arrow[r,hook]\arrow[u,hook] \arrow[uur, hook, bend right=4, crossing over] & 
    \I^*(X) \arrow[u,hook] \arrow[uur, hook, bend right=4, crossing over] \\
\end{tikzcd}\end{equation*}
Taking the double complexes $J^*$, $K^*$ corresponding to $\J^{*,*}$, $\K^{*,*}$, we recognize that the maps relating the cohomologies computed with $\A^*$ and $\B^*$ commute with the descent $\overline{\phi(X)}$ in cohomology:
\begin{equation*} \label{eq:AcyclicMorCoh} \begin{tikzcd} [row sep=0]
    \HH^k(\A^*(X)) \arrow[dd,"\overline{\phi(X)}","\sim"'{sloped}] \arrow[r,"\sim"] & \HH^k(J^*) \arrow[dd,"\overline{\phi(X)}","\sim"'{sloped}] & \\
    && \HH^k(\I^*(X))  \arrow[ul,"\sim"{sloped}] \arrow[dl,"\sim"{sloped}] \\
    \HH^k(\B^*(X)) \arrow[r,"\sim"] & \HH^k(K^*) &
\end{tikzcd}\end{equation*}
Thus, the isomorphism $\HH^k(\A^*(X))\xrightarrow{\sim} \HH^k(\B^*(X))$ induced by sheaf cohomology is equal to $\overline{\phi(X)}$.

\subsection{Refinement of \Cech{} cohomology} \label{CechRefinementQuasiIso}

An example of the quasi-isomorphism between acyclic resolutions that will interest us in \Cref{RealTori} is the refinement between \Cech{} resolutions.
Let $\U=\{U_i\}_i$ and $\V=\{V_j\}_j$ be two good open covers of $X$, with $\V$ a refinement of $\U$, meaning that there is a function $\rho:\V \to \U$ with $V\subseteq \rho(V)$ for all $V\in\V$. The restrictions $\F(\rho(V))\to \F(V)$ induce a chain map $\phi$ of the respective \Cech{} resolutions, which is the identity on $\F$ and a quasi-isomorphism  on global sections. So the induced map on cohomology $\oline{\phi(X)}$ is the isomorphism induced by the sheaf cohomology formalism.
\begin{equation*}
    \begin{tikzcd}[row sep=0]
        & \prod_{|\bmi|=*}\F_{|U_{\bmi}} \arrow[dd,"\phi"] &
        & \prod_{|\bmi|=*}\F(U_{\bmi}) \arrow[dd,"\phi(X)"] &[.4cm] 
        \cHH^*_{\U}(X,\F) \arrow[dd,"\sim"'{sloped},"\oline{\phi(X)}"] \\
        \F \arrow[ur,hook] \arrow[dr,hook] & & 
        \F(X) \arrow[ur,hook] \arrow[dr,hook] & & \\
        & \prod_{|\bmj|=*}\F_{|V_{\bmj}} &
        & \prod_{|\bmj|=*}\F(V_{\bmj}) &
        \cHH^*_{\V}(X,\F)
    \end{tikzcd}
\end{equation*}
In particular, it follows that the sheaf cohomology maps $\cHH^k_{\U}(X,\F) \xrightarrow{\sim} \HHH^k(X,\F)$ glue to a map from the direct limit over all good open covers 
$ \cHH^k(X,\F) \coloneq \displaystyle\lim_{\U\,\to} \cHH^k_{\U}(X,\F) \xrightarrow{\sim} \HHH^k(X,\F). $

\subsection{\Cech{} double complex}

In general, two acyclic resolutions need not come equipped with a natural chain map between them. In case one of them is the \Cech{} resolution, induced by a good cover $\U$ of $X$, one can find an intermediate acyclic resolution relating to both by considering the following diagrams, where the horizontal maps are the \Cech{} differentials $\cdelta$ and the vertical maps are the restrictions of the differentials $d_{\A}$ of $\A$:
\begin{equation} \label{eq:CechDiagram} \begin{tikzcd} [column sep=small, row sep=small]
    \A^* \arrow[r,hook] & \prod_{|\bmi|=*}\A^*_{|U_{\bmi}}
    &{}& 
    \A^*(X) \arrow[r,hook] & \prod_{|\bmi|=*}\A^*(U_{\bmi}) \\[.1cm]
    \F \arrow[r,hook]\arrow[u,hook] & \prod_{|\bmi|=*}\F_{|U_{\bmi}} \arrow[u,hook]
    &{}& 
    \F(X) \arrow[r,hook]\arrow[u,hook] & \prod_{|\bmi|=*}\F(U_{\bmi})  \arrow[u,hook] \\
\end{tikzcd}\end{equation}
The double complex $\E^*=\bigoplus_{a+b=*} \prod_{|\bmi|=a}\A^b_{|U_{\bmi}}$ with differential acting on $\prod_{|\bmi|=a}\A^b_{|U_{\bmi}}$ by $\cdelta + (-1)^a d_{\A}$ admits chain maps $\A^*\overset{\phi}{\hookrightarrow}\E^*\overset{\psi}{\hookleftarrow} \Cc_{\U}^*$ which induce quasi-isomorphisms $\A^*(X)\xrightarrow{\text{q.i.s.}}\E^*(X) \xleftarrow{\text{q.i.s.}} \Cc_{\U}^*(X)$. The above discussion then implies that the isomorphism $\HH^k(\A^*(X))\xrightarrow{\sim} \HH^k(\Cc_{\U}^*(X))$ induced by sheaf cohomology is given by $\overline{\psi(X)}^{-1}\circ \overline{\phi(X)}$. Concretely, this isomorphism is obtained by chasing cocycles across the anti-diagonals of diagram \eqref{eq:CechDiagram}.


\section{\large Proof of the theorem}

In \Cref{HomologicalSection} we have explained that the isomorphism $\HHs^*(X,\uZZ)\xrightarrow{\sim}\cHH_{\U}^*(X,\uZZ)$ induced by sheaf cohomology
is obtained by chasing cocycles across the anti-diagonals of the following commutative diagram:
\begin{equation} \label{eq:SinCohSquare} \begin{tikzcd}[]
    {} & {} & {} & {} & {} \\
    \tCs^{m+1}(X) \arrow[u,dotted,"\tdelta"] \arrow[r,"\cdelta"]  &  \PO_i \tCs^{m+1}(U_i) \arrow[u,dotted,"\tdelta"] \arrow[r,dotted,"\cdelta"]  &  \POO{k} \tCs^{m+1}(U_{\bmi}) \arrow[u,dotted,"\tdelta"] \arrow[r,"\cdelta"]  &  \POO{k+1} \tCs^{m+1}(U_{\bmi}) \arrow[u,dotted,"\tdelta"] \arrow[r,dotted,"\cdelta"]  & {} \\
    \tCs^{m}(X) \arrow[u,"\tdelta"] \arrow[r,"\cdelta"]  &  \PO_i \tCs^{m}(U_i) \arrow[u,"\tdelta"] \arrow[r,dotted,"\cdelta"]  &  \POO{k} \tCs^{m}(U_{\bmi}) \arrow[u,"\tdelta"] \arrow[r,"\cdelta"]  &  \POO{k+1} \tCs^{m}(U_{\bmi}) \arrow[u,"\tdelta"] \arrow[r,dotted,"\cdelta"]  & {} \\
    \tCs^{0}(X) \arrow[u,dotted,"\tdelta"] \arrow[r,"\cdelta"]  &  \PO_i \tCs^{0}(U_i) \arrow[u,dotted,"\tdelta"] \arrow[r,dotted,"\cdelta"]  &  \POO{k} \tCs^{0}(U_{\bmi}) \arrow[u,dotted,"\tdelta"] \arrow[r,"\cdelta"]  &  \POO{k+1} \tCs^{0}(U_{\bmi}) \arrow[u,dotted,"\tdelta"] \arrow[r,dotted,"\cdelta"]  & {} \\
    \uZZ(X) \arrow[u,"\tdelta"] \arrow[r,"\cdelta"]  &  \PO_i \uZZ(U_i) \arrow[u,"\tdelta"] \arrow[r,dotted,"\cdelta"]  &  \POO{k} \uZZ(U_{\bmi}) \arrow[u,"\tdelta"] \arrow[r,"\cdelta"]  &  \POO{k+1} \uZZ(U_{\bmi}) \arrow[u,"\tdelta"] \arrow[r,dotted,"\cdelta"]  & {} 
\end{tikzcd}\end{equation}

We would like to find an expression for the diagram chase. To this end, we want to reduce the singular cochains portion of diagram \eqref{eq:SinCohSquare} to the  $\Delta$-complex cohomology of the simplicial complex $N\U$. With the homotopy equivalence $\iota:N\U\to X$ we can do so only for the first column, as there are no meaningful restrictions of $\iota$ to homotopy equivalences $Y_{\bmi}\rightarrow U_{\bmi}$ for subcomplexes $Y_{\bmi}\subseteq N\U$. That is why we will reduce to the case where these restrictions exist, and are realized by enlarging the cover $\U$.

\begin{definition}
    For every open subset $O\subseteq X$, we define the restricted open cover $\U_{|O}=\{U\in\U \, | \,U\subseteq O\}$ and its associated \emph{subnerve} $N\U_{|O}\subseteq N\U$.
    For $O=U_{\bmi}$ a non-empty finite intersection in $\U$, we denote the subnerve shortly as  $N\U_{\bmi}$.
\end{definition}



\begin{definition}
    We call $\U$ \emph{saturated} if every non-empty finite intersection $U_{\bmi}$ is an open set already contained in $\U$.
\end{definition}

If $\U$ is saturated, every non-empty finite intersection $U_{\bmi}=U_{\bmi_0}\cap\cdots\cap U_{\bmi_{|\bmi |}}$ equals to $U_j\in\U$ for precisely one index $j$; for each $\bmi$ we denote this index by $\h{\bmi}$.
\footnote{It holds $U_{\bmi}=U_{\h{\bmi}}$ and likewise $\U_{\bmi}=\U_{\h{\bmi}}$ for the restricted covering; in particular $i=\h{\{i\}}$.}
For notational convenience, we assume without loss of generality that the ordering of the index set of $\U$ is such that $i<j$ for every $U_i\supset U_j$.
\footnote{This means that for any two $\bmj\subseteq\bmi$ it holds  $\h{\bmj}\leq\h{\bmi}$. One verifies that if $\bmj_0\leq\cdots\leq\bmj_k$ is ordered and $U_i\supseteq U_{\bmj_l},l=0,...,|\bmj|$, then $i\leq\bmj_0\leq\cdots\leq\bmj_k$ is also ordered. }
It's immediate that if $U_{\bmj}$ is a non-empty finite intersection with $U_{\bmj_l}\in\U_{\bmi}, l=0,...,|\bmj|$, then $U_{\bmi\cup\bmj}$ is also non-empty.
This translates to the fact that each subnerve $N\U_{\bmi}=N\U_{\h{\bmi}}\subset N\U$ has $x_{\h{\bmi}}\in N\U_{\bmi}$ as a so-called \emph{cone-vertex}, meaning that for every $k$-simplex $(x_{\bmj_0} \cdots x_{\bmj_{m}})\in\dC_k(N\U_{\bmi})$ with $\h{\bmi}\notin \bmj$ there is the $k+1$-simplex $(x_{\h{\bmi}}\,x_{\bmj_0} \cdots x_{\bmj_{m}})\in\dC_{k+1}(N\U_{\bmi})$.

\begin{lemma} \label[lemma]{SaturatedLemma}
    Assume $\U$ is saturated. 
    Then, for every non-empty finite intersection $U_{\bmi}$, the homotopy equivalence $\iota:N\U\rightarrow X$ restricts to a homotopy equivalence $N\U_{\bmi}\to U_{\bmi}$.
\end{lemma}
\begin{proof}
    Each subnerve $N\U_{\bmi}$ is non-empty and contractible to the cone-vertex $x_{\h{\bmi}}$. Since $U_{\bmi}$ is also contractible, $\iota:N\U_{\bmi}\to U_{\bmi}$ is a homotopy equivalence.
\end{proof}

We will reduce to the case where $\U$ is saturated using the following lemma, whose proof relies on the paracompactness hypothesis.
\begin{lemma} \label{CoverRefinement}
    Let $\U\subseteq \V$ be two ordered good open covers of $X$, with the ordering on $\V$ being an extension of the one on $\U$. Then, $N\U$ is naturally a homotopically equivalent subcomplex of $N\V$ and any homotopy equivalence $\iota:N\U\rightarrow X$ that respects $\U$ can be extended to a homotopy equivalence ${\iota'}:N\V\to X$ with ${\iota'}_{|N\U}=\iota$ that respects $\V$.
\end{lemma}
\begin{proof}
    The homotopy equivalence in \cite[Theorem 2.1.(1)]{paris} is constructed inductively on the skeleton $(N\U^n)_{n\geq 0}$ of the \Cech{} nerve. Thus, if we impose ${\iota'}_{|N\U^n}=\iota_{|N\U^n}$ at each step of the induction for $N\V$, we still get a homotopy equivalence ${\iota'}:N\V\to X$ that respects $\V$ and moreover extends $\iota$. 
    \[\begin{tikzcd}[row sep=normal]
        N\V \arrow[rrd,pos=0.3,"{\iota'}"] && \\
        N\U \arrow[u,hook] \arrow[rr,pos=0.4,"\iota"] && X
    \end{tikzcd} \]
\end{proof}

\subsection{Reduction to a bigger open cover. } \label{ReductionBiggerCover}
We want to justify that to prove \Cref{thm:Theorem} for the good cover $\U$ and $\iota$, it is sufficient to prove it for a bigger good cover $\V \supset \U$ and a homotopy $\iota':N\V\to X$ that extends $\iota$, viewing $N\U$ as a subcomplex of $N\V$.
The inclusion $\U\hookrightarrow\V$ induces a restriction homomorphism, which is a quasi-isomorphism of chain complexes, between diagram \eqref{eq:SinCohSquare} and the same diagram formed by the cover $\V$:

\begin{equation}\begin{tikzcd} []
    {} & {} & {} \\
    \tCs^{m+1}(X)\arrow[u,dotted,"\tdelta"] \arrow[rr, bend left=10,pos=0.1,dotted] \arrow[dr,pos=0.1,dotted,"\cdelta"]  & {} & 
    \DS\prod_{|\bmi|=*}\tCs^{m+1}(U_{\bmi}) \arrow[u,dotted,"\tdelta"] \\[-.1cm]
    \tCs^m(X) \arrow[u,"\tdelta"] \arrow[rr, bend left=20,pos=0.1,dotted] \arrow[dr,pos=0.1,dotted,"\cdelta"]   & 
    \DS\prod_{|\bmj|=*}\tCs^{m+1}(V_{\bmj}) \arrow[u,dotted,"\tdelta"] \arrow[ur, two heads,"\text{q.i.s.}"] &
    \DS\prod_{|\bmi|=*}\tCs^{m}(U_{\bmi}) \arrow[u,"\tdelta"]  \\
    \uZZ(X) \arrow[u,dotted,"\tdelta"] \arrow[rr, bend left=20,pos=0.1,dotted] \arrow[dr,pos=0.1,dotted,"\cdelta"]  &
    \DS\prod_{|\bmj|=*}\tCs^{m}(V_{\bmj}) \arrow[u,"\tdelta"] \arrow[ur, two heads,"\text{q.i.s.}"] & \check{C}_{\U}^*(X) \arrow[u,dotted,"\tdelta"]  \\
    & \check{C}_{\V}^*(X) \arrow[u,dotted,"\tdelta"] \arrow[ur, two heads,"\text{q.i.s.}"] & 
\end{tikzcd}\end{equation}

We recognize that chasing the sheafification $\tilde{\alpha}\in \tCs^m(X)$ of a singular cocycle $\alpha\in \Cs^m(X)$ across the anti-diagonals of diagram \eqref{eq:SinCohSquare} for $\V$, and then restricting to the \Cech{} complex of $\U$, results in a cocycle which is the diagram chase of $\tilde{\alpha}$ across the diagram for $\U$.
So if \Cref{thm:Theorem} holds for $\V$ and ${\iota'}$, after composing with the restriction we obtain
{\small$$ [\alpha] \in \HHs^k(X,\ZZ) \;\mapsto\;
    (-1)^{k}[\alpha({\iota'}\circ (x_{\bmj_0}\cdots x_{\bmj_k}))_{\bmj} ] \in\cHH_{\V}^k(X,\uZZ) \;\mapsto\;
    (-1)^{k}[\alpha(\iota\circ (x_{\bmi_0}\cdots x_{\bmi_k}))_{\bmi} ] \in\cHH_{\U}^k(X,\uZZ) ,  $$
\normalsize}
since ${\iota'}\circ (x_{\bmi_0}\cdots x_{\bmi_k}) = \iota\circ (x_{\bmi_0}\cdots x_{\bmi_k})$ for all index tuples $\bmi$ of $\U$.

\subsection{The case where $\U$ is saturated. } \label{CaseUSaturated}
With diagram \eqref{eq:SinCohSquare} we also have the ones obtained by replacing in it $\tCs^*(X)$ by $\Cs^*(X)$, $\Cs^*(N\U)$, $\Cd^*(N\U)$ and $\tCs^*(U_{\bmi})$ by $\Cs^*(U_{\bmi})$, $\Cs^*(N\U_{\bmi})$, $\Cd^*(N\U_{\bmi})$.

By \ref{ReductionBiggerCover} and \Cref{CoverRefinement}, we can assume without loss of generality that $\U$ is saturated.
Since every intersection $U_{\bmi}\subseteq X$ is likewise semi-locally contractible, every sheafification map $\Cs^*(U_{\bmi})\to \tCs^*(U_{\bmi})$ is a quasi-isomorphism \cite{sella}. 
By \Cref{SaturatedLemma}, the pullbacks $\iota_{|N\U_{\bmi}}^*:\Cs^*(U_{\bmi})\to \Cs^*(N\U_{\bmi})$ are quasi-isomorphisms, and so also the restriction to $\Delta$-complex cohomology.
Moreover,  for all $m\in\NN$, it is verified directly that the rows 
$\prod_{|\bmi|=*}\Cs^m(\U_{\bmi})$, $\prod_{|\bmi|=*}\Cs^m(N\U_{\bmi})$ and $\prod_{|\bmi|=*}\Cd^m(N\U_{\bmi})$ are exact.
So each of the four commutative diagrams has exact rows and columns, except for the first row and first column.

\begin{equation*} \label{eq:QuasiIso}
    \begin{tikzcd}[row sep = small, column sep = normal]
    & & & \tCs^*(X) \arrow[r,dotted,"\cdelta"]  &  \POO{*}\tCs^*(U_{\bmi}) \arrow[r,dotted, shorten >=.4cm, pos=0.3, "\cdelta"]  & {} \\
    & & \Cs^*(X) \arrow[r,dotted,"\cdelta"] \arrow[ur,"\text{sheafification}"] \arrow[dl,"\iota^*"']  &  \POO{*}\Cs^*(U_{\bmi})\arrow[r,dotted, shorten >=.4cm, pos=0.3, "\cdelta"]  \arrow[ur]  \arrow[dl]  & {} & \\ 
    & \Cs^*(N\U) \arrow[r,dotted,"\cdelta"] \arrow[dl,"\text{restriction}"']  &  \POO{*}\Cs^*(N\U_{\bmi})\arrow[r,dotted, shorten >=.4cm, pos=0.3, "\cdelta"]  \arrow[dl]   &  {} & & \\
    \Cd^*(N\U) \arrow[r,dotted,"\cdelta"]  &  \POO{*}\Cd^*(N\U_{\bmi})\arrow[r,dotted, shorten >=.4cm, pos=0.3, "\cdelta"]  &  {} & & & 
    \\[-1.6cm]
    & & & \uZZ(X) \arrow[uuuu,dotted,shorten >= 2.2cm,pos=0.2,"\tdelta"'] \arrow[r,dotted,"\cdelta"]  &  \POO{*} \uZZ(U_{\bmi}) \arrow[uuuu,dotted, shorten >= 1.6cm,pos=0.2,"\tdelta"'] \arrow[r,dotted, shorten >=.4cm, pos=0.3, "\cdelta"]  &  {} \\
    & & \uZZ(X) \arrow[dl,"\sim"{sloped}] \arrow[ur,"\sim"{sloped}] \arrow[uuuu,dotted, shorten >= 2.2cm,pos=0.2,"\delta"'] \arrow[r,dotted,"\cdelta"]  &  \POO{*} \uZZ(U_{\bmi}) \arrow[dl,"\sim"{sloped}] \arrow[ur,"\sim"{sloped}] \arrow[uuuu,dotted, shorten >= 1.6cm,pos=0.2,shift right=4,"\delta"']\arrow[r,dotted, shorten >=.4cm, pos=0.3, "\cdelta"]  & {} & \\
    & \uZZ(N\U) \arrow[dl,"\sim"{sloped}] \arrow[uuuu,dotted, shorten >= 2.2cm,pos=0.2,"\delta"'] \arrow[r,dotted,"\cdelta"]  &  \POO{*} \uZZ(N\U_{\bmi}) \arrow[dl,"\sim"{sloped}] \arrow[uuuu,dotted, shorten >= 1.6cm,pos=0.2,shift right=4,"\delta"']\arrow[r,dotted, shorten >=.4cm, pos=0.3, "\cdelta"]  & {} & & \\
    \uZZ(N\U) \arrow[uuuu,dotted, shorten >= 2.2cm,pos=0.2,"\delta"'] \arrow[r,dotted,"\cdelta"]  &  \POO{*} \uZZ(N\U_{\bmi}) \arrow[uuuu,dotted, shorten >= 1.6cm,pos=0.2,shift right=4,"\delta"']\arrow[r,dotted, shorten >=.4cm, pos=0.3, "\cdelta"]  & {} & & &
\end{tikzcd}\end{equation*}

With the facts mentioned in the previous paragraph, we recognize that chasing the sheafification $\tilde{\alpha}\in \tCs^m(X)$ of a singular cocycle $\alpha\in \Cs^m(X)$ across the anti-diagonals of diagram \eqref{eq:SinCohSquare} is equivalent, meaning up to coboundaties, to chasing $\iota^*(\alpha)\in \Cd^m(N\U)$ through the following diagram:

\begin{equation} \label{deltaCohSquare} \begin{tikzcd}
    {} & {} & {} & {} & {} \\
    \Cd^{m+1}(N\U) \arrow[u,dotted,"\delta"] \arrow[r,"\cdelta"]  &  \PO_i \Cd^{m+1}(N\U_i) \arrow[u,dotted,"\delta"] \arrow[r,dotted,"\cdelta"]  &  \POO{k} \Cd^{m+1}(N\U_{\bmi}) \arrow[u,dotted,"\delta"] \arrow[r,"\cdelta"]  &  \POO{k+1} \Cd^{m+1}(N\U_{\bmi}) \arrow[u,dotted,"\delta"] \arrow[r,dotted,"\cdelta"]  & {} \\
    \Cd^{m}(N\U) \arrow[u,"\delta"] \arrow[r,"\cdelta"]  &  \PO_i \Cd^{m}(N\U_i) \arrow[u,"\delta"] \arrow[r,dotted,"\cdelta"]  &  \POO{k} \Cd^{m}(N\U_{\bmi}) \arrow[u,"\delta"] \arrow[r,"\cdelta"]  &  \POO{k+1} \Cd^{m}(N\U_{\bmi}) \arrow[u,"\delta"] \arrow[r,dotted,"\cdelta"]  & {} \\
    \Cd^{0}(N\U) \arrow[u,dotted,"\delta"] \arrow[r,"\cdelta"]  &  \PO_i \Cd^{0}(N\U_i) \arrow[u,dotted,"\delta"] \arrow[r,dotted,"\cdelta"]  &  \POO{k} \Cd^{0}(N\U_{\bmi}) \arrow[u,dotted,"\delta"] \arrow[r,"\cdelta"]  &  \POO{k+1} \Cd^{0}(N\U_{\bmi}) \arrow[u,dotted,"\delta"] \arrow[r,dotted,"\cdelta"]  & {} \\
    \uZZ(N\U) \arrow[u,"\delta"] \arrow[r,"\cdelta"]  &  \PO_i \uZZ(N\U_i) \arrow[u,"\delta"] \arrow[r,dotted,"\cdelta"]  &  \POO{k} \uZZ(N\U_{\bmi}) \arrow[u,"\delta"] \arrow[r,"\cdelta"]  &  \POO{k+1} \uZZ(N\U_{\bmi}) \arrow[u,"\delta"] \arrow[r,dotted,"\cdelta"]  & {} 
\end{tikzcd}\end{equation}


\subsection{$\Delta$-complex and \Cech{} homology} \label{DeltaCechHomology}

Diagram \eqref{deltaCohSquare} is the dualization of
\begin{equation} \label{eq:DeltaHomSquare} \begin{tikzcd}[]
    {}\arrow[d,dotted,"\partial"] & {}\arrow[d,dotted,"\partial"]  & {} \arrow[d,dotted,"\partial"]  & {}\arrow[d,dotted,"\partial"] & {} \\
    \dC_{m+1}(N\U) \arrow[d,"\partial"]  &  \BO_i \dC_{m+1}(N\U_i) \arrow[d,"\partial"] \arrow[l,"\sum^{\Delta}"] 
    &  \BOO{k} \dC_{m+1}(N\U_{\bmi}) \arrow[d,"\partial"] \arrow[l,dotted,"\cpartial"] 
     &  \BOO{k+1} \dC_{m+1}(N\U_{\bmi}) \arrow[d,"\partial"] \arrow[l,"\cpartial"] 
      & {}\arrow[l,dotted,"\cpartial"]  \\
    \dC_{m}(N\U) \arrow[d,dotted,"\partial"]  &  \BO_i \dC_{m}(N\U_i) \arrow[d,dotted,"\partial"] \arrow[l,"\sum^{\Delta}"] 
     &  \BOO{k} \dC_{m}(N\U_{\bmi}) \arrow[d,dotted,"\partial"] \arrow[l,dotted,"\cpartial"] 
    &  \BOO{k+1} \dC_{m}(N\U_{\bmi}) \arrow[d,dotted,"\partial"] \arrow[l,"\cpartial"] 
    & {}\arrow[l,dotted,"\cpartial"] \\
    \dC_{0}(N\U) \arrow[d,"\sum_{N\U}"]  &  \BO_i \dC_{0}(N\U_i) \arrow[d,"(\sum_{N\U_i})_i"] \arrow[l,"\sum^{\Delta}"] 
    &  \BOO{k} \dC_{0}(N\U_{\bmi}) \arrow[d,"(\sum_{N\U_{\bmi}})_{\bmi}"] \arrow[l,dotted,"\cpartial"] 
    &  \BOO{k+1} \dC_{0}(N\U_{\bmi}) \arrow[d,"(\sum_{N\U_{\bmi}})_{\bmi}"] \arrow[l,"\cpartial"] 
    & {}\arrow[l,dotted,"\cpartial"] \\
    \uZZ(N\U)  &  \BO_i \uZZ(N\U_i) \arrow[l,"\sum"] 
    &  \BOO{k} \uZZ(N\U_{\bmi}) \arrow[l,dotted,"\cpartial"] 
    &  \BOO{k+1} \uZZ(N\U_{\bmi}) \arrow[l,"\cpartial"] 
    & {}\arrow[l,dotted,"\cpartial"] 
\end{tikzcd}\end{equation}

at least if we imagine the bottom row of \eqref{deltaCohSquare} being made up by the canonically isomorphic groups $\prod_{|\bmi|=*} \Hom(\uZZ(N\U_{\bmi}),\ZZ) $.
The vertical and horizontal maps of \eqref{eq:DeltaHomSquare} are the $\Delta$-complex boundary map and the covariant \Cech{} differential respectively.
Each complex $\dC_*(N\U_{\bmi})$ augmented by $\uZZ(N\U_{\bmi})$ in degree $-1$ is exact; note that the boundary operators $\partial$ in degree $0$ are simply the summation of the integral coefficients of the $0$-chains. Likewise, the rows except for the first one are exact and the covariant \Cech{} differential in degree $0$ is simply the summation of $0$-chains of the subcomplexes $N\U_i\subset N\U$.

We would like to find an expression for the diagram chase of a \Cech{} cycle in $\bigoplus_{|\bmi |=k}\uZZ(N\U_{\bmi})$ up to a simplicial cycle in $\dC_k(N\U)$ across the anti-diagonal, so we search for an operator $ \dC_*(N\U_{\bmi}) \to \dC_{*+1}(N\U_{\bmi}) $ which is a right-inverse of $\partial$ on cycles. With the crucial assumption that $\U$ is saturated, recall that each subnerve $N\U_{\bmi}$ has a cone-vertex $x_{\h{\bmi}}$, giving us a well-defined cone map
\begin{equation} \label{eq:ConeDef}
    \C_{x_{\h{\bmi}}}:  (x_{\bmj_0} \cdots x_{\bmj_{k}})\in\dC_k(N\U_{\bmi}) \mapsto
    \begin{cases} (x_{\h{\bmi}}\,x_{\bmj_0} \cdots x_{\bmj_{k}}) & \h{\bmi}\notin \bmj \\ 0 & \h{\bmi}\in\bmj \end{cases}
    \in \dC_{k+1}(N\U_{\bmi}),
\end{equation}
extended in degree $-1$ as $z\in\uZZ(N\U_{\bmi}) \mapsto z\cdot(x_{\h{\bmi}}) \in \dC_0(N\U_{\bmi})$. We denote $\C$ the generic cone operator on the columns, except the first one, of diagram \eqref{eq:DeltaHomSquare}. The next result is verified by an immediate computation.

\begin{lemma}
    $ \partial \C + \C \partial = \id, \quad \cdelta \,\C^{\vee} + \C^{\vee} \cdelta = \id $.
\end{lemma}
    


To compute the diagram chase across the anti-diagonals of diagram \eqref{eq:DeltaHomSquare}, we will need a closed expression for $(\cpartial \C)^{k+1}$.
The case distinction \eqref{eq:ConeDef} is handled formally, meaning that we allow ourselves to use the notation $(x_{v_0}\cdots x_{v_k})$, for a $k+1$-tuple of ordered indices $(v_0,...,v_k)$, even if some of these indices repeat, with the understanding that the expression $(x_{v_0}\cdots x_{v_k})$ will denote the $0$ element in $\dC_k(N\U)$. Note that, since $v$ is ordered, any repeating indices must be consecutive.
\begin{proposition} \label[proposition]{ithCone}
    For every generator $e_{\bmi}$ of the $\bmi$-th summand in $\bigoplus_{|\bmi|=k} \uZZ(N\U_{\bmi})$ holds
    \begin{multline} \label{eq:BarycentricSubd}
        (\cpartial \C)^{k+1}(\,e_{\bmi}\,) = \sum_{\sigma\in S\{0,...,k\}} \text{sgn}\sigma\,
            (x_{\h{\{\bmi_{\sigma (k)}\}}} \,
            x_{\h{\{\bmi_{\sigma (k)},\bmi_{\sigma (k-1)}\}}} \,
            x_{\h{\{\bmi_{\sigma (k)},\bmi_{\sigma (k-1)},\bmi_{\sigma (k-2)}\}}} \cdots \,
            x_{\h{\{\bmi_{\sigma (k)},...,\bmi_{\sigma (0)}\}}} ) \\
            = (-1)^{\frac{k(k+1)}{2}} \sum_{\sigma\in S\{0,...,k\}} \text{sgn}\sigma\,
            (x_{\h{\{\bmi_{\sigma (0)}\}}} \,
            x_{\h{\{\bmi_{\sigma (0)},\bmi_{\sigma (1)}\}}} \,
            x_{\h{\{\bmi_{\sigma (0)},\bmi_{\sigma (1)},\bmi_{\sigma (2)}\}}} \cdots \,
            x_{\h{\{\bmi_{\sigma (0)},...,\bmi_{\sigma (k)}\}}} ),
    \end{multline}
    where $(-1)^{\frac{k(k+1)}{2}}$ is the sign of the palyndromic permutation $i\mapsto k-i$ factored out.
\end{proposition}

\newcommand{\myvec}[2]{
    \begin{blockarray}{ c@{\hspace{0pt}} @{\hspace{8pt}} c}
        \begin{block}{ ( c@{\hspace{0pt}}) @{\hspace{8pt}} c}
          & \matindex{\vdots} \\
        #1 & \matindex{$#2$} \\
         & \matindex{\vdots} \\
        \end{block}
    \end{blockarray}
}
\newcommand{\myvecc}[4]{
    \begin{blockarray}{ c@{\hspace{0pt}} @{\hspace{8pt}} c}
        \begin{block}{ ( c@{\hspace{0pt}}) @{\hspace{8pt}} c}
          & \matindex{\vdots} \\
        #1 & \matindex{$#2$} \\
        #3 & \matindex{$#4$} \\
         & \matindex{\vdots} \\
        \end{block}
    \end{blockarray}
}
\newcommand{\myveccc}[6]{
    \begin{blockarray}{ c@{\hspace{0pt}} @{\hspace{8pt}} c}
        \begin{block}{ ( c@{\hspace{0pt}}) @{\hspace{8pt}} c}
          & \matindex{\vdots} \\
        #1 & \matindex{$#2$} \\
        #3 & \matindex{$#4$} \\
        #5 & \matindex{$#6$} \\
         & \matindex{\vdots} \\
        \end{block}
    \end{blockarray}
}

\begin{proof}
    First we illustrate the computation for $k=0,1,2$.

    \begin{minipage}{.3\textwidth}
        \small\[\begin{tikzcd}[row sep=normal, column sep=small]
            (x_i) &
            \myvec{(x_{i})}{\{i\}} \arrow[l,mapsto] \\
            & e_{\{i\}} \arrow[u,mapsto]
        \end{tikzcd}\] \normalsize
    \end{minipage}
    \begin{minipage}{.65\textwidth}
        \small\[\begin{tikzcd}[row sep=small, column sep=small]
        \begin{array}{c}
            -(x_{\h{\{\bmi_0\}}}\,x_{\h{\{\bmi_0,\bmi_1\}}}) \\
            +(x_{\h{\{\bmi_1\}}}\,x_{\h{\{\bmi_0,\bmi_1\}}})
        \end{array}
        & \myvecc{
            -(x_{\h{\{\bmi_0\}}}\,x_{\h{\{\bmi_0,\bmi_1\}}})}{\{\bmi_0\}}
            {(x_{\h{\{\bmi_1\}}}\,x_{\h{\{\bmi_0,\bmi_1\}}})}{\{\bmi_1\}} \arrow[l,mapsto] & \\
        & \myvecc{
            -(x_{\h{\{\bmi_0,\bmi_1\}}})}{\{\bmi_0\}}
            {(x_{\h{\{\bmi_0,\bmi_1\}}})}{\{\bmi_1\}} \arrow[u,mapsto]
            & \myvec{(x_{\h{\{\bmi_0,\bmi_1\}}})}{\{\bmi_0,\bmi_1\}} \arrow[l,mapsto] \\
            & & e_{\{\bmi_0,\bmi_1\}} \arrow[u,mapsto]
        \end{tikzcd}\] \normalsize
    \end{minipage}
                
    \tiny\[ \begin{tikzcd}[row sep=tiny, column sep=small,outer sep=-1pt]
        \begin{array}{c}
            (x_{\h{\{\bmi_0\}}}\, x_{\h{\{\bmi_0,\bmi_2\}}}\, x_{\h{\{\bmi_0,\bmi_1,\bmi_2\}}}) \\
            - (x_{\h{\{\bmi_0\}}}\, x_{\h{\{\bmi_0,\bmi_1\}}}\, x_{\h{\{\bmi_0,\bmi_1,\bmi_2\}}}) \\
            + (x_{\h{\{\bmi_1\}}}\, x_{\h{\{\bmi_0,\bmi_1\}}}\, x_{\h{\{\bmi_0,\bmi_1,\bmi_2\}}}) \\
            - (x_{\h{\{\bmi_1\}}}\, x_{\h{\{\bmi_1,\bmi_2\}}}\, x_{\h{\{\bmi_0,\bmi_1,\bmi_2\}}}) \\
            + (x_{\h{\{\bmi_2\}}}\, x_{\h{\{\bmi_1,\bmi_2\}}}\, x_{\h{\{\bmi_0,\bmi_1,\bmi_2\}}}) \\
            - (x_{\h{\{\bmi_2\}}}\, x_{\h{\{\bmi_0,\bmi_2\}}}\, x_{\h{\{\bmi_0,\bmi_1,\bmi_2\}}})
        \end{array}
        & \myveccc{
            (x_{\h{\{\bmi_0\}}}\, x_{\h{\{\bmi_0,\bmi_2\}}}\, x_{\h{\{\bmi_0,\bmi_1,\bmi_2\}}}) - (x_{\h{\{\bmi_0\}}}\, x_{\h{\{\bmi_0,\bmi_1\}}}\, x_{\h{\{\bmi_0,\bmi_1,\bmi_2\}}})}{\{\bmi_0\}}
            {(x_{\h{\{\bmi_1\}}}\, x_{\h{\{\bmi_0,\bmi_1\}}}\, x_{\h{\{\bmi_0,\bmi_1,\bmi_2\}}}) - (x_{\h{\{\bmi_1\}}}\, x_{\h{\{\bmi_1,\bmi_2\}}}\, x_{\h{\{\bmi_0,\bmi_1,\bmi_2\}}})}{\{\bmi_1\}}
            {(x_{\h{\{\bmi_2\}}}\, x_{\h{\{\bmi_1,\bmi_2\}}}\, x_{\h{\{\bmi_0,\bmi_1,\bmi_2\}}}) - (x_{\h{\{\bmi_2\}}}\, x_{\h{\{\bmi_0,\bmi_2\}}}\, x_{\h{\{\bmi_0,\bmi_1,\bmi_2\}}})}{\{\bmi_2\}}\arrow[l,mapsto]
        & & \\
        & \myveccc{
            (x_{\h{\{\bmi_0,\bmi_2\}}} \, x_{\h{\{\bmi_0,\bmi_1,\bmi_2\}}})-(x_{\h{\{\bmi_0,\bmi_1\}}} \, x_{\h{\{\bmi_0,\bmi_1,\bmi_2\}}})}{\{\bmi_0\}}
            {(x_{\h{\{\bmi_0,\bmi_1\}}} \, x_{\h{\{\bmi_0,\bmi_1,\bmi_2\}}})-(x_{\h{\{\bmi_1,\bmi_2\}}} \, x_{\h{\{\bmi_0,\bmi_1,\bmi_2\}}})}{\{\bmi_1\}}
            {(x_{\h{\{\bmi_1,\bmi_2\}}} \, x_{\h{\{\bmi_0,\bmi_1,\bmi_2\}}})-(x_{\h{\{\bmi_0,\bmi_2\}}} \, x_{\h{\{\bmi_0,\bmi_1,\bmi_2\}}})}{\{\bmi_2\}} \arrow[u,mapsto]
        & \myveccc{
            (x_{\h{\{\bmi_0,\bmi_1\}}} \, x_{\h{\{\bmi_0,\bmi_1,\bmi_2\}}})}{\{\bmi_0,\bmi_1\}}
            {-(x_{\h{\{\bmi_0,\bmi_2\}}} \, x_{\h{\{\bmi_0,\bmi_1,\bmi_2\}}})}{\{\bmi_0,\bmi_2\}}
            {(x_{\h{\{\bmi_1,\bmi_2\}}} \, x_{\h{\{\bmi_0,\bmi_1,\bmi_2\}}})}{\{\bmi_1,\bmi_2\}} \arrow[l,mapsto] & \\
        & \myveccc{
            (x_{\h{\{\bmi_0,\bmi_1,\bmi_2\}}})}{\{\bmi_0,\bmi_1\}}
            {-(x_{\h{\{\bmi_0,\bmi_1,\bmi_2\}}})}{\{\bmi_0,\bmi_2\}}
            {(x_{\h{\{\bmi_0,\bmi_1,\bmi_2\}}})}{\{\bmi_1,\bmi_2\}} \arrow[ur,mapsto]
        & \myvec{(x_{\h{\{\bmi_0,\bmi_1,\bmi_2\}}})}{\{\bmi_0,\bmi_1,\bmi_2\}} \arrow[l,mapsto] & \\
        & & e_{\{\bmi_0,\bmi_1,\bmi_2\}} \arrow[u,mapsto] &
    \end{tikzcd}\]
    \normalsize

    For the general case, we recognize that the result is a sum over all possible orders in which to remove indices from  $\{\bmi_{0},...,\bmi_{k}\}$ one by one. Removing indices in the stantard order $\bmi_{0}$, then $\bmi_{1}$, until $\bmi_{k-1}$ yields the term $(x_{\h{\{\bmi_{k}\}}} \,
            x_{\h{\{\bmi_{k},\bmi_{k-1}\}}} \,
            x_{\h{\{\bmi_{k},\bmi_{k-1},\bmi_{k-2}\}}} \cdots \,
            x_{\h{\{\bmi_{k},...,\bmi_{0}\}}} )$
    with sign $1$. Moreover, swapping two elements in the order of removal swaps also the sign, thus determining the sign of every other term.
\end{proof}

\begin{wrapfigure}{r}{0.34\textwidth}
  \begin{center}
    \includegraphics[width=0.33\textwidth]{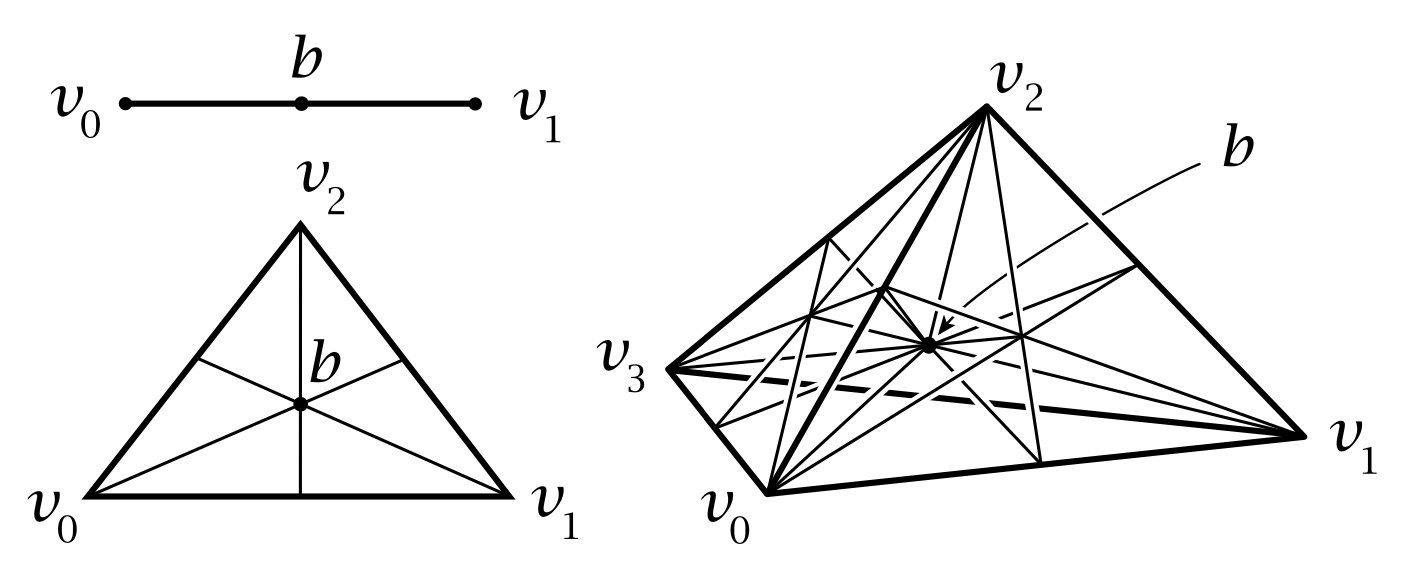}
    \cite[Proposition 2.21]{hatcher}
  \end{center}
\end{wrapfigure}
We recognize the term on the right-hand side of \eqref{eq:BarycentricSubd} to be a kind of barycentric subdivision of the $k$-simplex $(x_{\bmi_0}\cdots x_{\bmi_k})$, up to sign; to make this precise, we have to go back to the singular chains setting. Recall that the barycentric subdivision of a singular simplex $s\in\sC_k(X)$ is given by the signed sum of restrictions
$ S(s) := \sum_{\sigma\in S\{0,...,k\}} \text{sgn}\sigma\,s \circ \Delta^k_{\sigma}, $
where $\Delta^k_{\sigma}:\Delta^k\hookrightarrow\Delta^k$ is the embedded simplex of the standard $k$-simplex $\Delta^k=(v_0\cdots v_k)$ which has vertices $v_{\sigma0},v_{\h{\sigma0,\sigma1}},...,v_{\h{\sigma0,...,\sigma k}}$ in order, $v_{\h{i_0,...,i_j}}$ being the barycenter of the face $(v_{i_0}\cdots v_{i_j}) \subseteq \Delta^k$.
In \cite[Proposition 2.21]{hatcher} it is proved that there exists a homotopy operator $T:\sC_*(X)\to\sC_{*+1}(X)$ such that $T\partial+\partial T = id - S$.

In our $\Delta$-complex setting, the situation is analogous. For the operator
$$ S: (x_{\bmi_0}\cdots x_{\bmi_k}) \in\dC_k(N\U)\mapsto
    \sum_{\sigma\in S\{0,...,k\}} \text{sgn}\sigma\,
            (x_{\h{\{\bmi_{\sigma (0)}\}}} \,\cdots \,
            x_{\h{\{\bmi_{\sigma (0)},...,\bmi_{\sigma (k)}\}}} )
    \in\dC_k(N\U) $$
one can construct a homotopy operator $T:\dC_*(N\U)\to\dC_{*+1}(N\U)$ with $T\partial+\partial T = id - S$.

\subsection{Conclusion. }

Let $[\alpha]\in\HHs^k(X,\ZZ)$ and $[\iota^*(\alpha)]\in\HH_{\Delta}^k(N\U,\ZZ)$ its restriction to the cohomology of the simplicial complex $N\U$; see \eqref{eq:QuasiIso}.
By \ref{CaseUSaturated}, chasing the sheafification $\tilde{\alpha}\in \tCs^k(X,\uZZ)$ across diagram \eqref{eq:SinCohSquare} is equivalent to chasing $\iota^*(\alpha)$ across diagram \eqref{deltaCohSquare} and then across the isomorphism $\prod_{|\bmi|=k}\uZZ(N\U_{\bmi})\cong \prod_{|\bmi|=k}\uZZ(U_{\bmi})$.

By \ref{DeltaCechHomology}, the dual $\C^{\vee}$ is a cone operator for the columns of diagram \eqref{deltaCohSquare}.
Denote $\C'^{\,\vee}$ the cone operator acting on the column $\prod_{|\bmi |=a}\Cd^{*}(N\U_{\bmi})$ by $(-1)^a \C^{\vee}$.
Then, in the double complex $\prod_{|\bmi |=*}\Cd^{*}(N\U_{\bmi})$ with differential acting on $\prod_{|\bmi |=a}\Cd^{b}(N\U_{\bmi})$ by $\cdelta + (-1)^a \delta$, the cocycle 
$ \cdelta\iota^*(\alpha)\in \prod_i\Cd^k(N\U_i) $ is cohomologous to  
\footnote{Here we use $\sum_{0\leq i\leq k} i = \frac{k(k+1)}{2}$.}
\begin{equation} \label{6389}
    (-\, \cdelta\,\C'^{\,\vee})^{k} \,\cdelta \,\iota^*(\alpha) = 
    (-1)^{k+\frac{k(k+1)}{2}} \,\cdelta\, (\C^{\vee} \cdelta )^k \,\iota^*(\alpha)
    \in \prod_{|\bmi |=k}\Cd^{0}(N\U_{\bmi}).
\end{equation}

Evaluating at the cone vertices $\C(e_{\bmi})\in C^{\Delta}_0(N\U_{\bmi})$ and using \Cref{ithCone}, for $e_{\bmi}$ the generator of the $\bmi$-th summand in $\bigoplus_{|\bmi|=k}\uZZ(N\U_{\bmi})$, we obtain the constant cocycle in $\prod_{|\bmi|=k} \uZZ(N\U_{\bmi})$ corresponding to \eqref{6389}:
\begin{equation}\label{eq:ChasedCocycle} 
    (-1)^{k+\frac{k(k+1)}{2}}\left( \iota^*(\alpha) (\cpartial \C)^{k+1}(e_{\bmi})  \right)_{\bmi}  
    = (-1)^{k} \left( \alpha(\iota \circ S(x_{\bmi_0}\cdots x_{\bmi_k})) \right)_{\bmi}
    \in \prod_{|\bmi|=k} \uZZ(N\U_{\bmi}),
\end{equation}
where $S$ is the barycentric subdivision operator.
Let's use the homotopy operator $T$ described in \ref{DeltaCechHomology} to relate the cocycle \eqref{eq:ChasedCocycle} to the one claimed in the theorem:
\begin{multline*}
    (\alpha(\iota\circ (x_{\bmi_0} \cdots x_{\bmi_{k}}))
        - \alpha(\iota\circ S(x_{\bmi_0} \cdots x_{\bmi_{k}})) )_{\bmi}
    = (\alpha(\iota\circ T\partial (x_{\bmi_0} \cdots x_{\bmi_{k}})) )_{\bmi}
        + \cancel{ (\alpha(\iota\circ \partial T (x_{\bmi_0} \cdots x_{\bmi_{k}})) )_{\bmi}}
    \\
    = \left(\sum_{l=0}^{k}(-1)^l \alpha(\iota\circ T(x_{\bmi_0}\cdots \h{x_{\bmi_l}} \cdots x_{\bmi_{k}})) \right)_{\bmi}
    = \cdelta \left(\,(\alpha(\iota\circ T(x_{\bmj_0}\cdots x_{\bmj_{k-1}})))_{\bmj}\, \right)
    \; \in \POOj{k}\uZZ(N\U_{\bmj}).
\end{multline*}

Consequently, the two \Cech{} cocycles are cohomologous. So $(-1)^{k}(\alpha(\iota\circ (x_{\bmi_0} \cdots x_{\bmi_{k}})))_{\bmi}$ likewise represents the chasing of $\tilde{\alpha}$ across diagram \eqref{eq:SinCohSquare} 
, i.e. a representative of the \Cech{} cohomology class corresponding to $[\alpha]$. The theorem is thus proven.

\section{\large Automorphy and \Cech{} Maps of Real Tori} \label[section]{RealTori}

\newcommand{\lambdaSin}{\lambda^{\mathrm{sin}}}
\newcommand{\lambdaSinCo}{\lambda_{\mathrm{sin}}}
\newcommand{\lambdaDelta}{\lambda^{\Delta}}
\newcommand{\lambdaDeltaCo}{\lambda_{\Delta}}
\newcommand{\lambdaNabla}{\lambda^{\nabla}}
\newcommand{\lambdavee}{\lambda^{\vee}}

\newcommand{\olinelambdavee}{\oline{\lambda}^{\vee}}

\newcommand{\latticevec}[1]{\lambda_{#1}}
\newcommand{\latticedual}[1]{\lambda^{#1}}
\newcommand{\latticedualbar}[1]{\oline{\lambda}^{#1}}

After establishing notation for singular, \Cech{} and group cohomology on a real torus, we use Theorem \ref{thm:Theorem} to show that the automorphy map $\HHs^k(V/\Lambda, \ZZ) \xrightarrow{\sim} \HH^k(\Lambda, \HH^0(V,\pi^*\uZZ))$ and the group cohomology map $\HH^k(\Lambda, \HH^0(V,\pi^*\uZZ)) \xrightarrow{\sim} \cHH^k_{\U}(V/\Lambda,\uZZ)$ induced by the universal cover  $\pi:V \twoheadrightarrow V/\Lambda$ are inverse to one another, if we choose the good cover $\U$ in a sensible way with respect to the lattice of the torus.

\subsection{Cohomology of real tori}
Let $\Lambda=\left<\latticevec{i},\dots,\latticevec{n}\right>$ be a free abelian group of rank $n$, $V=\Lambda\otimes\RR$ the $\RR$-vector space containing $\Lambda$ as lattice and $X=V/\Lambda$ the $n$-dimensional torus. We denote by 
$$ \latticedual{i} :  \sum_j a_j\latticevec{j} \in V \mapsto a_i \in\RR,
\qquad \latticedualbar{i}: \sum_j \oline{a_j}\latticevec{j} \in X \mapsto a_i \in [0,1) $$
the linear coordinate functions; $\latticedual{i}$ will also be considered as a linear form in $\HomZ{\Lambda}$ or $\Hom_{\RR}(V,\RR)$.
The torus $X$ has (co)homology bases in degree $1$ consisting of 
\begin{align*}
     \Lambda & \xrightarrow{\sim} \sHH_1(X,\ZZ), &\; \latticevec{i} & \mapsto [\lambdaSin_i], &\; \left( \lambdaSin_i: s\in[0,1]\mapsto \oline{s}\in X \right) \;\; \in \sCC_1(X,\ZZ) , 
\end{align*}
\begin{align*}
    \Hom(\Lambda,\ZZ) & \xrightarrow{\sim} \HHs^1(X,\ZZ), &\; \latticedual{i} & \mapsto [\lambdaSinCo^i], &\; \lambdaSinCo^i
    \in \CCs^1(X,\ZZ) & \;\text{ with }\; \lambdaSinCo^i(\lambdaSin_j) = \delta_{ij} \text{ for all $j$}.
\end{align*}
We have the Künneth isomorphisms, where $\Alt^k(\Lambda,\ZZ)$ the group of $k$-alternating forms on $\Lambda$:
\begin{equation*} \label{eq:KunnethDiagram} \begin{tikzcd}[row sep=large]
    \Alt^k(\Lambda,\ZZ)  \dar[hook] & \lar["\sim"']
    \bigwedge^k\Hom(\Lambda,\ZZ) \rar["\wedge_j \latticedual{\bmi_j} \,\mapsto\, \wedge_j {[\lambdaSinCo^{\bmi_j}]}"{yshift=0.2cm},"\sim"'] \dar[hook] &[.2cm]
    \bigwedge^k\HHs^1(X,\ZZ) \rar["\sim"',"\smile"] &
    \HHs^k(X,\ZZ) \arrow[dr,hook,xshift=.3cm,yshift=.15cm] &
    \\
    \Alt^k(\Lambda,\RR) &  \lar["\sim"']
    \bigwedge^k\Hom(\Lambda,\RR) \rar["\wedge_j \latticedual{\bmi_j} \,\mapsto\, \wedge_j {[d\latticedual{\bmi_j}]}"{yshift=.2cm},"\sim"'] &
    \bigwedge^k\HHdR^1(X,\uRR) \rar["\sim"',"\wedge"] &
    \HHdR^k(X,\uRR) \rar["\sim"',"{[\omega]} \,\mapsto\, {[\int_{\bullet}\omega]}"{yshift=.1cm,xshift=-.2}] &
    \HHs^k(X,\RR)
\end{tikzcd}\end{equation*}

\subsection{Automorphy map}

The constant automorphy factors of primitive functions of closed differential forms on $V$ with respect to the deck transformation group $\mathrm{Deck}(\pi)\cong\pi_1(X)\cong\Lambda$ are equal 
to the periods of differential forms on $X$. This computation can be carried out at the level of singular cohomology, giving an isomorphism 
\begin{equation*}\begin{tikzcd}[column sep=huge]
    \HHs^1(X,\ZZ)
    \arrow[rr, bend left=5,"\sim"'{yshift=-.05cm},"{[\alpha]} \,\mapsto\, {[\lambda \,\mapsto\, \pi^*\alpha \,( 0\to\lambda ) ]} "]  & &
    \HHH^1(\Lambda,\ZZ) \arrow[ll, bend left=5, "{[\lambdaSinCo^i] \mapsfrom {\lambda^i}}"] \,,
\end{tikzcd}\end{equation*}
where $(0\to\lambda):s\in[0,1]\mapsto s\cdot \lambda \in V$ are any singular $1$-simplices on $V$ starting in $0$ and ending in $\lambda$.
A generalization of this map, which we will call the \emph{automorphy} map, is obtained through the corresponding cup products in singular and group cohomology:
\begin{equation}\begin{tikzcd}[column sep=huge] \label{AutomorphyDiagram}
    \HHs^k(X,\ZZ)
    \arrow[rr,"\sim"',"{[\alpha]} \,\mapsto\, {[(l_1,\dots,l_k) \,\mapsto\, \pi^*\alpha \,(0\to l_1\to l_1+l_2\to\dots\to\sum_j l_j ) ]} "] & &
    \HHH^k(\Lambda,\ZZ)
    \\
    \bigwedge^k \HHs^1(X,\ZZ)
    \arrow[rr,"\sim"',"\wedge_j\,{[\alpha_j]} \,\mapsto\, \wedge_j\,{[\lambda\mapsto \pi^*\alpha_j\,( 0 \to \lambda ) ]} "]
    \uar["\sim"{sloped},"\smile"'] & &
    \bigwedge^k \HHH^1(\Lambda,\ZZ)
    \uar["\sim"{sloped},"\hspace{-.08cm}\mapstoUP{\wedge_j\,{[F_{\bmi_j}]}}{ {[ (l_1,\dots,l_k) \,\mapsto\, F_{\bmi_1}(l_1)\cdots F_{\bmi_k}(l_k) ]} }"']
    \\
    & \bigwedge^k \HomZ{\Lambda}
    \arrow[ul,"\sim"'{sloped}]
    \arrow[ur,"\sim"{sloped},"\wedge_j\, \latticedual{\bmi_j} \,\mapsto\, \wedge_j\, {[\latticedual{\bmi_j}]} "'] &
\end{tikzcd}\end{equation}
where $(0\to l_1\to l_1+l_2\to\dots\to\sum_j l_j )$ are any singular $k$-simplices on $V$ with vertices in order $0$, $l_1$, $l_1+l_2$, \dots, $\sum_j l_j$.
The inverse isomorphism from group cohomology to alternating forms is given by alternization:
\begin{equation*}\begin{tikzcd}[column sep=normal]
    \Alt^k(\Lambda,\ZZ) \rar["\sim"] &
    \bigwedge^k\HomZ{\Lambda} 
    \arrow[rr,"\substack{ \wedge_j \latticedual{\bmi_j} \,\mapsto\, \wedge_j{[\latticedual{\bmi_j}]} \\ \sim}"{yshift=.08cm}]  & &
    \bigwedge^k \HHH^1(\Lambda,\ZZ) \arrow[rr,"\substack{ \wedge_j{[F_{\bmi_j}]} \,\mapsto\, {[ (l_1,\dots,l_k) \,\mapsto\, F_{\bmi_1}(l_1)\cdots F_{\bmi_k}(l_k) ]} \\ \sim}"{yshift=0.08cm}] &[1.6cm] &
    \HHH^k(\Lambda,\ZZ) \,, \arrow[lllll,bend left=7,shift left=1,"\sim"',"\sum_{\sigma\in S(\{1,\dots,k\})} \sgn(\sigma)\,F\circ P_{\sigma} \;\mapsfrom\, {[F]}"{yshift=-.2cm}]
\end{tikzcd}\end{equation*}
where $P_{\sigma}(v_1 \dots,v_k)=(v_{\sigma 1},\dots,v_{\sigma k})$ is the permutation of the arguments. 

\subsection{\Cech{} and group cohomology}

\Cech{} cohomology of a sheaf of abelian groups $\F$ on $X$ is related via the universal cover $\pi: V\twoheadrightarrow X$ to group cohomology in the following way.
Let $\U=\{U_i\}_i$ be an ordered good open cover of $X$ and choose a connected component $\tilde{U_i}\subset \pi^{-1}(U_i)$ for all $i$, so that the restrictions $\pi_i:\tilde{U_i}\xrightarrow{\sim}U_i$ are homeomorphisms and, for all $U_i\cap U_j\neq \emptyset$, the differences $ \lambda_{ij} \coloneq \pi_j^{-1}-\pi_i^{-1} \in \Lambda$ are constant lattice elements; we are using the convention adopted in \cite[§2 Appendix]{mumford2008abelian} and \cite[§B]{BiLa}.
Then, for each $k$ there is an isomorphism
$$ \phi_{\U}^k : \HHH^k(\Lambda,\Gamma(V,\pi^*\F)) \xrightarrow{\sim} \cHH_{\U}^k(X,\F), \quad [F]\mapsto [(F(\lambda_{\bmi_0\bmi_1},\dots,\lambda_{\bmi_{k-1}\bmi_k}))_{\bmi}], $$
which is canonical in the sense that it is indipendent of the choices of components $\tilde{U_i}$.
\footnote{Different choices of components $\tilde{U_i},\tilde{U_i}'\subset \pi^{-1}(U_i)$ would result in lattice elements $\lambda_{ij},\lambda_{ij}'$ which differ by a \Cech{} coboundary $\lambda_{ij}'-\lambda_{ij} = ((\pi_j^{-1})'-\pi_j^{-1})-((\pi_i^{-1})'-\pi_i^{-1})$.}

Let now $\V=\{V_l\}_l$ be a refined good open cover of $\U$, meaning that there is a function $\rho:\V \to \U$ with $V\subseteq \rho(V)$ for all $V\in\V$.
Denoting $\rho l$ the index of $\rho(V_l)$ in $\U$, we can choose connected components $\tilde{V}_l\subset \pi^{-1}(V_l)$ such that $\tilde{V}_l\subseteq \tilde{U}_{\rho l}$.
Then, denoting $\pi_l:\tilde{V}_l\xrightarrow{\sim} V_l$ the restrictions, it holds $\pi^{-1}_m-\pi^{-1}_l = \pi^{-1}_{\rho m}-\pi^{-1}_{\rho l} = \lambda_{\rho l,\rho m} $.
So the refinement map on cohomology commutes with $\phi^k_{\U}$ and $\phi^k_{\V}$, and there is an induced isomorphism $\phi^k$ to the direct limit:  
\begin{equation*}
    \begin{tikzcd}[column sep=0cm,row sep=normal]
        \cHH_{\U}^k(X,\F) \arrow[rr,"\text{refinement}","\sim"'] && \cHH_{\V}^k(X,\F) \arrow[r,"\sim"] &[1.0cm] 
        \displaystyle\lim_{\U \,\to } \cHH_{\U}^k(X,\F) \arrow[r,phantom,"\eqcolon"] &[.5cm] \cHH^k(X,\F) \\
        & \HHH^k(\Lambda,\Gamma(V,\pi^*\F)) \arrow[ul,"\sim"{sloped,pos=0.5},pos=0.7,"\phi^k_{\U}"] \arrow[ur,"\sim"{sloped,pos=0.5},pos=0.7,"\phi^k_{\V}"'] \arrow[urr,bend right=8,"\sim"{sloped},"\phi^k"'] & \\
    \end{tikzcd}
\end{equation*}




\subsection{Simplicial complex structure}

We construct a simplicial complex structure on the torus $X$; the $2$-dimensional case is illustrated in diagram \eqref{N2Simplex}. 
As $0$-simplices we take the $3$-division points $(x_i)_{0\leq i< 3^n}$ of the torus, ordered lexicographically by the rule
$$ i<j :\iff \latticedualbar{k}(x_i)<\latticedualbar{k}(x_j), \text{ for } k=\text{max}\{1\leq l\leq n \,|\, \latticedualbar{l}(x_i)\neq \latticedualbar{l}(x_j)\} . $$
As $1$-simplices, for all $i<j$ such that 
$$ \left(\forall 1\leq k\leq n: 0\leq \latticedual{k}(x_j)-\latticedual{k}(x_i)\leq\frac{1}{3} \right) \text{ or } \left(\forall 1\leq k\leq n: -\frac{1}{3}\leq \latticedual{k}(x_j)-\latticedual{k}(x_i)\leq 0 \right) ,$$
we take the oriented segment $(x_i\to x_j)$ from $x_i$ to $x_j$ that stays within the $3$-division subcube of $X$ containing both vertices; so either $(x_i\to x_j)$ or its mirrored segment points positively in all directions.
The rest of the skeleton is uniquely determined by adding all possible higher dimensional simplices, which are necessarily each contained within one subcube. We denote a $k$-simplex $(x_{\bmi_0}\to\cdots\to x_{\bmi_k})$ by its ordered sequence of $0$-simplices.
Consider the following (dual to each other) bases for $\Delta$-complex (co)homology.

\begin{minipage}{.25\textwidth}
    \begin{align*}
        \Lambda & \xrightarrow{\sim} \dHH_1(X,\ZZ) \\
        \latticevec{k} & \mapsto [\lambdaDelta_k]
    \end{align*}
\end{minipage}
\begin{minipage}{.72\textwidth}
    $$ \lambdaDelta_k \coloneq \left(0 \to \frac{\latticevec{k}}{3}\right) + \left(\frac{\latticevec{k}}{3} \to \frac{2\latticevec{k}}{3}\right) - \left(0 \to \frac{2\latticevec{k}}{3}\right)  \in C^{\Delta}_1(X,\ZZ) $$
\end{minipage}

\begin{minipage}{.25\textwidth}
    \begin{align*}
        \HomZ{\Lambda} & \xrightarrow{\sim} \HHd^1(X,\ZZ)\\
        \latticedual{k} & \mapsto [\lambdaDeltaCo^k]
    \end{align*}
\end{minipage}
\begin{minipage}{.72\textwidth}
    $$ \lambdaDeltaCo^k : (x_i \to x_j)\in C^{\Delta}_1(X,\ZZ) \mapsto \begin{cases} 1 & \text{if } \latticedualbar{k}(x_i)=0,\,\latticedualbar{k}(x_j)=\frac{1}{3} \\
        -1 & \text{if } \latticedualbar{k}(x_i)=\frac{1}{3},\,\latticedualbar{k}(x_j)=0 \\
    0 & \text{else} \end{cases} $$
\end{minipage}

The described simplicial complex structure is realized in a natural way as the \Cech{} nerve $N\U$ of a good open cover $\U=\{U_i\}_{0\leq i< 3^n}$ of $X$ with 
$$ (x_{\bmi_0}\to\cdots\to x_{\bmi_k}) \subset \bigcup_{0\leq j\leq k} U_{\bmi_j}, \quad \text{ for all $k$-simplices}, $$
and ordered in the same way as $(x_i)_{0\leq i< 3^n}$.
We can choose connected components $\tilde{U}_i\subset \pi^{-1}(U_i)$ in such a way that 
the differences between the restrictions $\pi_i:\tilde{U}_i\xrightarrow{\sim}U_i$ are 
\begin{align} \label{LambdaIJequations}
    \Lambda\ni \lambda_{ij} = \pi_j^{-1}-\pi_i^{-1} = \sum_k
    \begin{cases}
        -\latticevec{k} & \text{ if } \latticedualbar{k}(x_i)=0,\,\latticedualbar{k}(x_j)=\frac{1}{3} \\
        \latticevec{k} & \text{ if } \latticedualbar{k}(x_i)=\frac{1}{3},\,\latticedualbar{k}(x_j)=0 \\
        0 & \text{ else}
    \end{cases} \;, \quad \text{for all } i<j,
\end{align}
letting us compute explicitely
$$ \latticedual{l}(\lambda_{ij}) = \begin{cases}
        -1 & \text{ if } \latticedualbar{l}(x_i)=0,\,\latticedualbar{l}(x_j)=\frac{1}{3} \\
        1 & \text{ if } \latticedualbar{l}(x_i)=\frac{1}{3},\,\latticedualbar{l}(x_j)=0 \\
        0 & \text{ else}
    \end{cases} \;, \quad \text{for all } l,\,i<j. $$
So $\latticedual{l}(\lambda_{ij}) = -\lambdaDeltaCo^l(x_i \to x_j)$ means that the following diagram commutes:
\begin{equation*} \label{eq:1dimDeltaDiagram}\begin{tikzcd}[row sep=large,column sep=huge]
    \HHd^1(X,\ZZ)    &&
    \HHH^1(\Lambda,\ZZ)
    \arrow[ll,"\sim"," {[(x_{i} \to x_{j})\mapsto F(\lambda_{ij})]} \,\mapsfrom \,\, {[F]} "'] 
    \\
    \HomZ{\Lambda}
    \uar["\sim"{sloped},"\mapstoUP{ \latticedual{i} }{ {[\lambdaDeltaCo^{i}]} }"']
    \arrow[rr,"\cdot (-1)","\sim"']
    && \HomZ{\Lambda}
    \uar["\sim"{sloped},"\mapstoUP{ \latticedual{i} }{ {[\latticedual{i}]} }"']
\end{tikzcd}\end{equation*}

Passing to the exterior product and applying Theorem \ref{thm:Theorem}, we obtain the following commutative diagram: 
\begin{equation*} \label{eq:CohDiagram} \begin{tikzcd}[row sep=huge,column sep=huge]
    \HHs^k(X,\uZZ) \rar["\substack{\text{sheaf}\\\text{cohomology}}","\sim"'] &[-.4cm]
    \cHH_{\U}^k(X,\uZZ) & \hspace{1.0cm} &
    \HHH^k(\Lambda,\ZZ)
    \arrow[ll,"\sim","{[F(\lambda_{\bmi_0\bmi_1},\dots,\lambda_{\bmi_{k-1}\bmi_k})_{\bmi} ]} \mapsfrom \,{[F]} \,:\, \phi^k_{\U}"']
    \\
    \HHs^k(X,\ZZ) \rar["\text{restriction}","\sim"'] \uar["\sim"{sloped},"\text{sheafification}"'] &
    \HHd^k(X,\ZZ)
    \uar["\sim"{sloped},"\mapstoUPleft{ \;\;\;[\alpha]}{ {[(\alpha(x_{\bmi_0}\to\cdots\to x_{\bmi_k}))_{\bmi}]}\cdot(-1)^{k} }"'{xshift=-.1cm}] &&
    \HHH^k(\Lambda,\ZZ)
    \arrow[ll,"\sim","{[(x_{\bmi_0}\to\cdots\to x_{\bmi_k})\mapsto F(\lambda_{\bmi_0\bmi_1},\dots,\lambda_{\bmi_{k-1}\bmi_k}) ]} \,\mapsfrom \,\,{[F]}"']
    \uar["\sim"{sloped},"\cdot(-1)^{k}"']
    \\
    \bigwedge^k \HHs^1(X,\ZZ) \rar["\text{restriction}","\sim"']
    \uar["\sim"{sloped},"\smile"'] &
    \bigwedge^k \HHd^1(X,\ZZ)
    \uar["\sim"{sloped},"\smile"']
    &&
    \bigwedge^k \HHH^1(\Lambda,\ZZ)
    \uar["\sim"{sloped},"\hspace{-.08cm}\mapstoUP{\wedge_j\,{[F_{\bmi_j}]}}{ {[ (l_1,\dots,l_k) \,\mapsto\, F_{\bmi_1}(l_1)\cdots F_{\bmi_k}(l_k) ]} }"']
    \arrow[ll,"\sim"," \wedge_j\,{[(x_{\mu}\to x_{\nu})\mapsto F_{\bmi_j}(\lambda_{\mu\nu})]} \,\mapsfrom \,\, \wedge_j\,{[F_{\bmi_j}]} "'] 
    \\
    & \bigwedge^k \HomZ{\Lambda}
    \arrow[lu,"\sim"'{sloped}," \wedge_j\, {[ \lambdaSinCo^{\bmi_j} ]} \,\mapsfrom\, \wedge_j\, \latticedual{\bmi_j} "{yshift=-.26cm,xshift=.3cm}]
    \uar["\sim"{sloped},"\mapstoUP{ \wedge_j\, \latticedual{\bmi_j} }{ \wedge_j\, {[\lambdaDeltaCo^{\bmi_j}]} }"']
    \arrow[rr,"\cdot (-1)^k","\sim"']
    && \bigwedge^k \HomZ{\Lambda}
    \uar["\sim"{sloped},"\mapstoUP{ \wedge_j\, \latticedual{\bmi_j} }{ \wedge_j\, {[\latticedual{\bmi_j}]} }"']
\end{tikzcd}\end{equation*}

De Rham and $\Delta$-complex cohomology fit in a similar commutative diagram:
\begin{equation*} 
\begin{tikzcd}[row sep=large, column sep=huge]
    \HHs^k(X,\uRR) \arrow[r,"\sim"'{sloped}] 
    & \cHH_{\U}^k(X,\uRR) & &[.45cm] 
    \HHH^k(\Lambda,\RR) \arrow[ll,"\sim"]
    \arrow[ll,"\sim","{[F(\lambda_{\bmi_0\bmi_1},\dots,\lambda_{\bmi_{k-1}\bmi_k})_{\bmi} ]} \mapsfrom \,{[F]} \,:\, \phi^k_{\U}"']  \\
    \HHdR^k(X,\uRR) \arrow[ru, bend left=0,"\substack{\text{sheaf}\\\text{cohomology}}"{xshift=-.1cm,yshift=-.1cm},"\sim"'{sloped}]
    \arrow[u,"\sim"'{sloped},"\mapstoUP{ {[\omega]} }{ {[\int_{\bullet}\omega ]} }"]
    \arrow[r,"\sim"'{sloped}," {[\omega]} \,\mapsto\, {[\int_{\bullet}\omega ]}"{xshift=.2cm}] 
    & \HHd^k(X,\RR)
    \uar["\sim"{sloped},"\mapstoUPleft{ \;\;\;[\alpha]}{ {[(\alpha(x_{\bmi_0}\to\cdots\to  x_{\bmi_k}))_{\bmi}]}\cdot(-1)^{k} }"'{xshift=-.1cm}] & &
    \\
    \bigwedge^k \HHdR^1(X,\uRR)
    \uar["\sim"{sloped},"\wedge"']
    \rar["\sim"',"\wedge_j\,{[\omega_j]} \,\mapsto\,\wedge_j\,{[\int_{\bullet}\omega_j ]} "{yshift=.05cm}]
    &
    \bigwedge^k \HHd^1(X,\RR)
    \uar["\sim"{sloped},"\smile"']
    & &
    \bigwedge^k \HHH^1(\Lambda,\RR)
    \arrow[uu,"\sim"{sloped},"\hspace{-.05cm}\mapstoUP{\wedge_j\,{[F_{\bmi_j}]}}{ {[ (l_1,\dots,l_k) \,\mapsto\, F_{\bmi_1}(l_1)\cdots F_{\bmi_k}(l_k) ]} }"']
    \\
    & \bigwedge^k \Hom(\Lambda,\RR)
    \arrow[lu,"\sim"'{sloped}," \wedge_j\, {[ d \latticedual{\bmi_j} ]} \,\mapsfrom\, \wedge_j\, \latticedual{\bmi_j} "{yshift=-.26cm,xshift=.3cm}]
    \uar["\sim"{sloped},"\mapstoUP{ \wedge_j\, \latticedual{\bmi_j} }{ \wedge_j\, {[\lambdaDeltaCo^{\bmi_j}]} }"']
    \arrow[urr,"\sim"{sloped},"\wedge_j\, \latticedual{\bmi_j} \,\mapsto\, \wedge_j\, {[\latticedual{\bmi_j}]} "'] & &
\end{tikzcd}\end{equation*}

Recalling that the automorphy map factors through $\bigwedge^k\Hom(\Lambda,\ZZ)$ as in diagram \eqref{AutomorphyDiagram}, we obtain that the following diagram commutes:
\begin{equation*}
    \begin{tikzcd}[column sep=normal]
        \HHs^k(V/\Lambda,\uZZ)
        \arrow[r,"\sim"',"\substack{\text{sheaf}\\\text{cohomology}}"] &[.8cm]
        \cHH^k_{\U}(V/\Lambda,\uZZ) \arrow[r,"\sim"',"\displaystyle\lim_{\U\,\to}"] &
        \cHH^k(V/\Lambda,\uZZ)
        \\[.2cm]
        \HHs^k(V/\Lambda,\ZZ) \arrow[u,"\sim"{sloped},"\text{sheafification}"'] \arrow[r,"\sim"'{sloped}," \text{automorphy}"] & 
        \HHH^k(\Lambda,\ZZ)
        \arrow[u,"\sim"{sloped},"\phi^k_{\U}"'] \arrow[ur,"\sim"{sloped},"\phi^k"'] 
    \end{tikzcd}
\end{equation*}

\section{\large Collating Formula}

Let us briefly recall the collating formula \cite[Proposition 9.8]{botttu} for a smooth manifold $X$.
If $(\rho_i)_i$ is a partition of unity on $X$ subordinate to $\U$, 
then there is a homotopy operator 
$$ K\,:\; (\omega_{\bmi})_{\bmi} \in \prod_{|\bmi|=k}\Omega_X^*(U_{\bmi}) \,\mapsto\, 
\left( \sum_{j} \sgn_{j,\bmj} \rho_j \,\omega_{\{j,\bmj_0,\dots, \bmj_{k-1}\}} \right)_{\bmj}  \in\prod_{|\bmj|=k-1}\Omega_X^*(U_{\bmj}) , \quad \cdelta\,K + K\cdelta = \id, $$
where $\sgn_{j,\bmj}$ is the sign of the permutation that orders $(j,\bmj_0,\dots, \bmj_{k-1})$.
The collating formula says that,
in the double complex $\prod_{|\bmi |=*}\Omega_X^{*}(U_{\bmi})$ with differential acting on $\prod_{|\bmi |=a}\Omega_X^{b}(U_{\bmi})$ by $\cdelta + (-1)^a d$, 
a constant \Cech{} cocycle $ (z_{\bmi})_{\bmi} \in \prod_{|\bmi|=k} \Omega_X^0(U_{\bmi}) $ is cohomologous to  
$(-1)^{k+\frac{k(k+1)}{2}} (d\, K )^{k} ((z_{\bmi})_{\bmi}) \in  \prod_i\Omega_X^k(U_i)$, which glues to a global $k$-form.

Our description of the sheaf cohomology isomorphism gives an inverse to this construction, or, in a sense, an explanation of it:  the de Rham cohomology class $[\omega]$ corresponding to a \Cech{} class $[(z_{\bmi})_{\bmi}]$ is the one represented by a global $k$-form $\omega$ which satisfies
\begin{equation} \label{1829}
    (-1)^{\frac{k(k+1)}{2}}\,\int_{ \iota\circ s_{\bmj}} (d\, K )^{k} ((z_{\bmi})_{\bmi}) = 
    (-1)^{k}\,\int_{ \iota\circ s_{\bmj}} \omega = 
    z_{\bmj},\quad \text{for all $k$-simplices $s_{\bmj}$ of $N\U$,} 
\end{equation}
under a homotopy equivalence $\iota:N\U\to X$ that respects $\U$.


\[\begin{tikzcd}
    \HHs^k(X,\uRR)
    \arrow[rr,"\sim"',"\substack{\text{sheaf}\\\text{cohomology}}"] &&
    \cHH^k_{\U}(X,\uRR) \arrow[ddr, bend left=34, "\sim"'{sloped}, "\substack{\text{collating } \text{formula} \\ (-1)^{k+\frac{k(k+1)}{2}} (d\, K )^{k} }"] & \\
    \HHs^k(X,\RR) \arrow[rr,"\iota^*","\sim"']\arrow[u,"\sim"{sloped},"\text{sheafification}"'] &&
    \HH_{\Delta}^k(N\U,\RR) \arrow[u,"\sim"{sloped},"\cdot (-1)^{k}"'] & \\
    & & & \HH_{dR}^k(X,\uRR)
    \arrow[lllu,bend left=12,"\sim"'{sloped},"\int_{\bullet}\omega\,\mapsfrom\,\omega"']
\end{tikzcd}\]




\subsection{Example} \label{TorusExample}

\usetikzlibrary{calc}
\usetikzlibrary{decorations,decorations.pathmorphing,decorations.pathreplacing,decorations.markings}
\tikzset{->-/.style={
    decoration={
        markings,
        mark=at position .55 with {\arrow{>}}},postaction={decorate}}}
\tikzstyle{every node}=[scale=1.2]

Let us apply \Cref{thm:Theorem} and verify \eqref{1829} in the case of a $2$-dimensional torus.
We illustrate the lattice $\langle\lambda_1,\lambda_2\rangle\subset V$, the simplicial structure of $X$ and the connected components $\tilde{U_i}\subset \pi^{-1}(\pi(\tilde{U_i}))$:
\begin{equation} \label{N2Simplex}
    \begin{tikzpicture}[x=1.5cm,y=1.5cm]
        \node[circle,fill=black,label=below right:{$0$}, inner sep=1pt] (l0) at (-5,0) {};
        \node[circle,fill=black,label=below right:{$\frac{\lambda_1}{3}$}, inner sep=1pt] (l1) at (-4,0) {};
        \node[circle,fill=black,label=below right:{$\frac{2\lambda_1}{3}$}, inner sep=1pt] (l2) at (-3,0) {};
        \node[circle,fill=black,label=below right:{$\lambda_1$}, inner sep=1pt] (l0') at (-2,0) {};
        \node[circle,fill=black,label=below right:{$\frac{\lambda_2}{3}$}, inner sep=1pt] (l3) at (-5,1) {};
        \node[circle,fill=black,label=below right:{}, label=above right:{\textcolor{red}{$\tilde{U}_4$}}, inner sep=1pt] (l4) at (-4,1) {};
        \node[circle,fill=black,label=below right:{}, label=above right:{\textcolor{orange}{$\tilde{U}_5$}}, inner sep=1pt] (l5) at (-3,1) {};
        \node[circle,fill=black,label=below right:{}, label=above right:{\textcolor{brown}{$\tilde{U}_3$}}, inner sep=1pt] (l3') at (-2,1) {};
        \node[circle,fill=black,label=below right:{$\frac{2\lambda_2}{3}$}, inner sep=1pt] (l6) at (-5,2) {};
        \node[circle,fill=black,label=below right:{}, label=above right:{\textcolor{blue}{$\tilde{U}_7$}}, inner sep=1pt] (l7) at (-4,2) {};
        \node[circle,fill=black,label=below right:{}, label=above right:{\textcolor{cyan}{$\tilde{U}_8$}}, inner sep=1pt] (l8) at (-3,2) {};
        \node[circle,fill=black,label=below right:{}, label=above right:{\textcolor{magenta}{$\tilde{U}_6$}}, inner sep=1pt] (l6') at (-2,2) {};
        \node[circle,fill=black,label=below right:{$\lambda_2$}, inner sep=1pt] (l0'') at (-5,3) {};
        \node[circle,fill=black,label=above right:{\textcolor{purple}{$\tilde{U}_1$}}, inner sep=1pt] (l1'') at (-4,3) {};
        \node[circle,fill=black,label=above right:{\textcolor{green}{$\tilde{U}_2$}}, inner sep=1pt] (l2'') at (-3,3) {};
        \node[circle,fill=black,label=below right:{$\lambda_1+\lambda_2$}, label=above right:{\textcolor{teal}{$\tilde{U}_0$}}, inner sep=1pt] (l0''') at (-2,3) {};
        \fill[red, fill opacity=0.2, rotate around={45:(l4)}] ($(l4)+(-1.4cm,-0.8cm)$) rectangle ($(l4)+(1.4cm,0.8cm)$);
        \fill[orange, fill opacity=0.2, rotate around={45:(l5)}] ($(l5)+(-1.4cm,-0.8cm)$) rectangle ($(l5)+(1.4cm,0.8cm)$);
        \fill[brown, fill opacity=0.2, rotate around={45:(l3')}] ($(l3')+(-1.4cm,-0.8cm)$) rectangle ($(l3')+(1.4cm,0.8cm)$);
        \fill[blue, fill opacity=0.2, rotate around={45:(l7)}] ($(l7)+(-1.4cm,-0.8cm)$) rectangle ($(l7)+(1.4cm,0.8cm)$);
        \fill[cyan, fill opacity=0.2, rotate around={45:(l8)}] ($(l8)+(-1.4cm,-0.8cm)$) rectangle ($(l8)+(1.4cm,0.8cm)$);
        \fill[pink, fill opacity=0.2, rotate around={45:(l6')}] ($(l6')+(-1.4cm,-0.8cm)$) rectangle ($(l6')+(1.4cm,0.8cm)$);
        \fill[purple, fill opacity=0.2, rotate around={45:(l1'')}] ($(l1'')+(-1.4cm,-0.8cm)$) rectangle ($(l1'')+(1.4cm,0.8cm)$);
        \fill[green, fill opacity=0.2, rotate around={45:(l2'')}] ($(l2'')+(-1.4cm,-0.8cm)$) rectangle ($(l2'')+(1.4cm,0.8cm)$);
        \fill[teal, fill opacity=0.2, rotate around={45:(l0''')}] ($(l0''')+(-1.4cm,-0.8cm)$) rectangle ($(l0''')+(1.4cm,0.8cm)$);
        %
        \node[circle,fill=black,label=below right:{$x_0$}, inner sep=1pt] (x0) at (0,0) {};
        \node[circle,fill=black,label=below right:{$x_1$}, inner sep=1pt] (x1) at (1,0) {};
        \node[circle,fill=black,label=below right:{$x_2$}, inner sep=1pt] (x2) at (2,0) {};
        \node (x0') at (3,0) {};
        \node[circle,fill=black,label=below right:{$x_3$}, inner sep=1pt] (x3) at (0,1) {};
        \node[circle,fill=black,label=below right:{$x_4$}, inner sep=1pt] (x4) at (1,1) {};
        \node[circle,fill=black,label=below right:{$x_5$}, inner sep=1pt] (x5) at (2,1) {};
        \node (x3') at (3,1) {};
        \node[circle,fill=black,label=below right:{$x_6$}, inner sep=1pt] (x6) at (0,2) {};
        \node[circle,fill=black,label=below right:{$x_7$}, inner sep=1pt] (x7) at (1,2) {};
        \node[circle,fill=black,label=below right:{$x_8$}, inner sep=1pt] (x8) at (2,2) {};
        \node (x6') at (3,2) {};
        \node (x0'') at (0,3) {};
        \node (x1'') at (1,3) {};
        \node (x2'') at (2,3) {};
        \node (x0''') at (3,3) {};
        \draw (x0.center) edge[->-] (x1.center); \draw (x1.center) edge[->-] (x2.center); \draw (x0'.center) edge[->-] (x2.center); 
        \draw (x3.center) edge[->-] (x4.center); \draw (x4.center) edge[->-] (x5.center); \draw (x3'.center) edge[->-] (x5.center); 
        \draw (x6.center) edge[->-] (x7.center); \draw (x7.center) edge[->-] (x8.center); \draw (x6'.center) edge[->-] (x8.center); 
        \draw (x0''.center) edge[->-] (x1''.center); \draw (x1''.center) edge[->-] (x2''.center); \draw (x0'''.center) edge[->-] (x2''.center); 
        \draw (x0.center) edge[->-] (x3.center); \draw (x3.center) edge[->-] (x6.center); \draw (x0''.center) edge[->-] (x6.center);
        \draw (x1.center) edge[->-] (x4.center); \draw (x4.center) edge[->-] (x7.center); \draw (x1''.center) edge[->-] (x7.center);
        \draw (x2.center) edge[->-] (x5.center); \draw (x5.center) edge[->-] (x8.center); \draw (x2''.center) edge[->-] (x8.center);
        \draw (x0'.center) edge[->-] (x3'.center); \draw (x3'.center) edge[->-] (x6'.center); \draw (x0'''.center) edge[->-] (x6'.center);
        \draw (x0.center) edge[->-] (x4.center); \draw (x1.center) edge[->-] (x5.center); \draw (x2.center) edge[->-] (x3'.center);
        \draw (x3.center) edge[->-] (x7.center); \draw (x4.center) edge[->-] (x8.center); \draw (x5.center) edge[->-] (x6'.center);
        \draw (x1''.center) edge[->-] (x6.center); \draw (x2''.center) edge[->-] (x7.center); \draw (x0'''.center) edge[->-] (x8.center);
    \end{tikzpicture}
\end{equation}

In concordance with equations \eqref{LambdaIJequations}, the only non-zero differences $\lambda_{ij}=\pi^{-1}_j-\pi^{-1}_i$ are: 
\begin{align*}
    \lambda_{01} &= -\lambda_1 \quad &
    \lambda_{03} &= -\lambda_2 \quad &
    \lambda_{04} &= -\lambda_1-\lambda_2 \\
    \lambda_{14} &= -\lambda_2 \quad&   
    \lambda_{15} &= -\lambda_2 \quad&
    \lambda_{16} &= \lambda_1 \\
    \lambda_{23} &= -\lambda_2 \quad&  
    \lambda_{25} &= -\lambda_2 \\
    \lambda_{34} &= -\lambda_1 \quad& 
    \lambda_{37} &= -\lambda_1 \\
    \lambda_{67} &= -\lambda_1
\end{align*}

Diagram \eqref{AutomorphyChechDiagram} says that the \Cech{} class corresponding to $\squareBr{\lambdaSinCo^i}\in\HHs^1(X,\ZZ)$ under the sheaf cohomology isomorphism is  $\phi^1_{\U}(\lambda^i)\in\cHH_{\U}^1(X,\ZZ)$, whose components are
\begin{equation*}
    \left(\phi^1_{\U}(\lambda^1)\right)_{ij}
     = \lambda^1(\lambda_{ij}) =
     \begin{cases}
        -1 & (i,j)=(0,1) \\
        -1 & (i,j)=(0,4) \\
        1 & (i,j)=(1,6) \\
        -1 & (i,j)=(3,4) \\
        -1 & (i,j)=(3,7) \\
        -1 & (i,j)=(6,7) \\
        0 & \text{ else} \\
     \end{cases}, \qquad
     \left(\phi^1_{\U}(\lambda^2)\right)_{ij}
     = \lambda^2(\lambda_{ij}) =
     \begin{cases}
        -1 & (i,j)=(0,3) \\
        -1 & (i,j)=(0,4) \\
        -1 & (i,j)=(1,4) \\
        -1 & (i,j)=(1,5) \\
        -1 & (i,j)=(2,3) \\
        -1 & (i,j)=(2,5) \\
        0 & \text{ else} \\
     \end{cases}\,.
\end{equation*}

Similarly, the \Cech{} class corresponding to the orientation form  
$\squareBr{\lambdaSinCo^1 \smile \lambdaSinCo^2}\in\HHs^2(X,\ZZ)$ is
\begin{equation*} \label{4895}
    \left(\phi^2_{\U}(\lambda^1\smile \lambda^2)\right)_{ijk}
    = (\lambda^1\smile \lambda^2)(\lambda_{ij},\lambda_{jk})
    = \lambda^1(\lambda_{ij})\cdot \lambda^2(\lambda_{jk}) = 
    \begin{cases}
    1 & (i,j,k)=(0,1,4) \\
    0 & \text{ else} \\
    \end{cases} \; \in\cHH_{\U}^2(X,\uZZ).
\end{equation*}

Then, we compute $(-1)^{2+\frac{2\cdot3}{2}}\cdot(d\,K)^2(\phi^2_{\U}(\lambda^1\smile \lambda^2)) = 2 d\rho_1\wedge d\rho_4 $, which is indeed a global form agreeing with $-2d\rho_0\wedge d\rho_4$ and $2d\rho_0\wedge d\rho_1$ on $U_{0,1,4}$,  so it is supported inside of the triple intersection $U_{0,1,4}$.

{\tiny\[ \begin{tikzcd}[row sep=-.1cm, column sep=normal]
    -2 d\rho_1\wedge d\rho_4 \arrow[r,mapsto]
    & \myveccc{-2 d\rho_1\wedge d\rho_4}{\{0\}}
        {2 d\rho_0\wedge d\rho_4}{\{1\}}
        {-2 d\rho_0\wedge d\rho_1}{\{4\}}
    & & \\ 
    & \myveccc{
        -\rho_1\cdot d\rho_4 + \rho_4\cdot d\rho_1}{\{0\}}
        {\rho_0\cdot d\rho_4 - \rho_4\cdot d\rho_0}{\{1\}}
        {\rho_1\cdot d\rho_0 - \rho_0\cdot d\rho_1}{\{4\}} \arrow[u,phantom, "\rotatebox{90}{$\mapsto$}\quad","d"']
    & \myveccc{
        d\rho_4 }{\{0,1\}}
        {-d\rho_1}{\{0,4\}}
        {d\rho_0}{\{1,4\}} \arrow[l,mapsto,"K"'] & \\
    & & \myveccc{
        \rho_4 }{\{0,1\}}
        {-\rho_1 }{\{0,4\}}
        { \rho_0 }{\{ 1,4\}} \arrow[u,phantom, "\rotatebox{90}{$\mapsto$}\quad","d"']
    & \myvec{1}{\{0,1,4\}} \arrow[l,mapsto,"K"']
\end{tikzcd}\]
\normalsize}
Looking back at our illustration \eqref{N2Simplex}, the good cover can be constructed such that a lift $\pi:\h{U} \xrightarrow{\sim} U_{0,1,4}$ of the triple intersection looks like a small square, for some factor $\epsilon\in\RR_{> 0}$:
\begin{equation*} \label{TripleIntersection}
    \begin{tikzpicture}[x=1.1cm,y=1.1cm]
        \node[circle,fill=black,label=right:{\small$v_0$},inner sep=1pt] (a) at (0,0) {};
        \node (b) at (-1,1) {};
        \node (c) at (1,1) {};
        \node (d) at (0,2) {};
        \node (f) at (0,1) {\tiny$\h{U}$};
        \draw (a.center) edge[->-]
            node[midway, below left] {\tiny$\epsilon(-\lambda_1+\lambda_2)$} (b.center);
        \draw (a.center) edge[->-]node[midway, below right] {\tiny$\epsilon(\lambda_1+\lambda_2)$}  (c.center);
        \draw (b.center) edge[->-]node[midway, above left] {\tiny$\epsilon(\lambda_1+\lambda_2)$}  (d.center); 
        \draw (c.center) edge[->-]node[midway, above right] {\tiny$\epsilon(-\lambda_1+\lambda_2)$}  (d.center); 
    \end{tikzpicture}
\end{equation*}
and the partition of unity can be chosen so that, for all $v\in\h{U}$, it satisfies 
$$ \rho_1(\pi(v - v_0)) = 1-\frac{-\lambda^1(v-v_0)+\lambda^2(v-v_0)}{2\epsilon}, \quad 
\rho_4(\pi(v - v_0)) = \begin{cases*}
    \frac{\lambda^2(v-v_0)}{\epsilon}-1 & \text{if }$\lambda^2(v-v_0)\geq \epsilon$, \\
    0 & \text{ else.}
\end{cases*} $$
Then, we compute 
$$ \int_{X} 2d\rho_1\wedge d\rho_4 = 
2\int_{ \left\{ v\in \h{U} | \lambda^2(v-v_0)\geq \epsilon\right\} } -\frac{-d\lambda^1 + d\lambda^2}{2\epsilon}\wedge \frac{d\lambda^2}{\epsilon}  = 
\frac{1}{\epsilon^2} \int_{  \left\{ v\in \h{U} | \lambda^2(v-v_0)\geq \epsilon\right\} } d\lambda^1\wedge d\lambda^2 
 = 1, $$
so $(-1)^{2+\frac{2\cdot3}{2}}\cdot(d\,K)^2(\phi^2_{\U}(\lambda^1\smile \lambda^2))$ is indeed the de Rham cohomology class corresponding to $\squareBr{\lambdaSinCo^1\smile \lambdaSinCo^2}$ because their integrals coincide.
Moreover, the equations \eqref{1829} are satisfied:
$$ (-1)^{\frac{2\cdot3}{2}} \int_{(x_i\to x_j\to x_k)} (d\,K)^2(\phi^2_{\U}(\lambda^1\smile \lambda^2)) = \begin{cases*}
    1 & \text{if $(i,j,k)=(0,1,4)$} \\
    0 & \text{else}
\end{cases*} \;
= \left( \phi^2_{\U}(\lambda^1\smile \lambda^2) \right)_{ijk}\,. $$

\bibliographystyle{alpha}
\bibliography{References.bib}

\end{document}